\documentclass[oneside, 12pt]{amsart}
\usepackage{amscd, amssymb, amsmath, mathrsfs}
\usepackage[english]{babel}
\usepackage{booktabs}
\usepackage{tikz-cd}
\usepackage[all]{xy}
\usepackage[pdftex, colorlinks=true,  citecolor=blue, linkcolor=blue, linktocpage=false]{hyperref}

\newcommand\mylabel[1]{\label{#1}\marginpar{\vspace{-1ex}\medskip\medskip\footnotesize \tt #1}}
\renewcommand\mylabel[1]{\label{#1}}
\newcommand{\mydate}{
\number\day\space
\ifcase\month \or January\or February\or March\or April\or May\or June\or July\or August\or September\or October\or November\or December\fi 
\space\number\year}

\DeclareUrlCommand\arXiv{\urlstyle{same}}

\newtheorem{theorem}{Theorem}[section]
\newtheorem{maintheorem}{Theorem}
\newtheorem{lemma}[theorem]{Lemma}
\newtheorem{proposition}[theorem]{Proposition}
\newtheorem{corollary}[theorem]{Corollary}

\theoremstyle{definition}

\newtheorem*{acknowledgement}{Acknowledgement}

\theoremstyle{remark}

\newcommand{\ZZ}{\mathbb{Z}}
\newcommand{\QQ}{\mathbb{Q}}
\newcommand{\RR}{\mathbb{R}}
\newcommand{\CC}{\mathbb{C}}
\newcommand{\FF}{\mathbb{F}}

\newcommand{\PP}{\mathbb{P}}

\newcommand{\GG}{\mathbb{G}}

\newcommand{\ideala}{\mathfrak{a}}

\newcommand{\idealg}{\mathfrak{g}}

\newcommand{\shA}{\mathscr{A}}

\newcommand{\shE}{\mathscr{E}}
\newcommand{\shF}{\mathscr{F}}

\newcommand{\shH}{\mathscr{H}}
\newcommand{\shHom}{\mathscr{H}\!\text{\textit{om}}}
\newcommand{\shI}{\mathscr{I}}

\newcommand{\shK}{\mathscr{K}}
\newcommand{\shM}{\mathscr{M}}

\newcommand{\shN}{\mathscr{N}}
\newcommand{\shL}{\mathscr{L}}
\newcommand{\shP}{\mathscr{P}}

\newcommand{\catC}{\mathcal{C}}

\newcommand{\foA}{\mathfrak{A}}

\newcommand{\foY}{\mathfrak{Y}}
\newcommand{\foZ}{\mathfrak{Z}}

\newcommand{\alg}{\text{\rm alg}}

\newcommand{\Alb}{\operatorname{Alb}}

\newcommand{\Art}{\text{\rm Art}}

\newcommand{\Aut}{\operatorname{Aut}}

\newcommand{\can}{\text{\rm  can}}

\newcommand{\CH}{\operatorname{CH}}

\newcommand{\ch}{\operatorname{ch}}

\newcommand{\Ext}{\operatorname{Ext}}

\newcommand{\Frac}{\operatorname{Frac}}

\newcommand{\GL}{\operatorname{GL}}

\newcommand{\gr}{\operatorname{gr}}

\newcommand{\Grp}{{\text{{\rm Grp}}}}

\newcommand{\Hom}{\operatorname{Hom}}
\newcommand{\Hilb}{\operatorname{Hilb}}

\newcommand{\kod}{\operatorname{Kod}}

\newcommand{\length}{\operatorname{length}}
\newcommand{\Lie}{\operatorname{Lie}}

\newcommand{\invlim}{\varprojlim}

\newcommand{\lra}{\longrightarrow}

\newcommand{\maxid}{\mathfrak{m}}

\newcommand{\NS}{\operatorname{NS}}

\renewcommand{\O}{\mathscr{O}}

\newcommand{\op}{\text{\rm op}}

\newcommand{\Pic}{\operatorname{Pic}}

\newcommand{\pr}{\operatorname{pr}}
\newcommand{\Proj}{\operatorname{Proj}}

\newcommand{\quadand}{\quad\text{and}\quad}

\newcommand{\ra}{\rightarrow}

\newcommand{\red}{{\operatorname{red}}}

\newcommand{\Spec}{\operatorname{Spec}}
\newcommand{\Spf}{\operatorname{Spf}}

\newcommand{\Sym}{\operatorname{Sym}}

\newcommand{\td}{\operatorname{td}}

\newcommand{\trdeg}{\operatorname{trdeg}}

\newcommand{\uHom}{\underline{\operatorname{Hom}}}

\newcommand  {\Lif}     {\operatorname{Lif}}

\begin{document}

\title[Moret-Bailly families]
      {Moret-Bailly families   and non-liftable schemes}

\author[Damian R\"ossler]{Damian R\"ossler}
\address{Mathematical Institute, University of Oxford, Andrew Wiles Building, Rad-
cliffe Observatory Quarter, Woodstock Road, Oxford, OX2 6GG, United Kingdom}
\curraddr{}
\email{rossler@maths.ox.ac.uk}

\author[Stefan Schr\"oer]{Stefan Schr\"oer}
\address{Mathematisches Institut, Heinrich-Heine-Universit\"at,
40204 D\"usseldorf, Germany}
\curraddr{}
\email{schroeer@math.uni-duesseldorf.de}

\subjclass[2010]{14J32, 14K05, 14L15, 14D10}
% 14L15 Group schemes
% 14D10 Arithmetic ground fields (finite, local, global) and families or fibrations
% 14J32 Calabi-Yau manifolds (algebro-geometric aspects) 
% 14K05 Algebraic theory of abelian varieties

\dedicatory{Revised Version, 29 June 2021}

\begin{abstract}
Generalizing the Moret-Bailly pencil of supersingular abelian surfaces to higher dimensions,
we  construct for each  field of characteristic $p>0$ a smooth projective variety  with trivial
dualizing sheaf that does not   lift to characteristic zero. 
Our approach heavily  relies on local unipotent group schemes,  the Beauville--Bogomolov Decomposition   
for K\"ahler manifolds with $c_1=0$, and equivariant deformation theory
in mixed characteristics.
\end{abstract}

\maketitle
\tableofcontents

\newcommand{\dashra}{\dashrightarrow}
\newcommand{\lieg}{\mathfrak{g}}
\newcommand{\PS}{{\PP^n}}
\newcommand{\frob}{{(p)}}
\newcommand{\uExt}{\underline{\operatorname{Ext}}}
\newcommand{\liea}{\mathfrak{a}}

%===========================================================
\section*{Introduction}
\mylabel{Introduction}

Every compact K\"ahler manifold $V$
with Chern class $c_1=0$  has unobstructed deformations,
although the obstruction group $H^2(V,\Theta_{V})$ is usually non-zero.
This foundational fact relies on the \emph{$T^1$-Lifting Theorem}
(confer \cite{Bogomolov 1978},  \cite{Todorov 1980}, \cite{Tian 1987},
\cite{Kawamata 1992},  \cite{Ran 1992}).
It holds, in particular, for complex tori, hyperk\"ahler manifolds, and
Calabi--Yau manifolds.  In fact,  the \emph{Beauville--Bogomolov Decomposition Theorem}
asserts that a compact K\"ahler manifold $V$ with $c_1=0$ admits a finite \'etale covering
$V'\ra V$ that splits into a product $V'=V_1\times\ldots\times V_r$ where the factors belong to these
three classes (\cite{Bogomolov 1974}, \cite{Beauville 1983}).

Over fields $k$ of characteristic $p>0$, much less is known for smooth proper scheme $Y$  
that have $c_1=0$, in the sense that the dualizing sheaf  $\omega_Y$ is numerically trivial.
Again we have three particular classes: The \emph{abelian varieties} take over the role of complex tori.
Copying the definition in characteristic zero, one may call  $Y$  a \emph{hyperk\"ahler manifold} if it is simply-connected, 
and there is some $\sigma\in H^0(Y,\Omega^2_Y)$ whose adjoint map  $\sigma:\Theta_Y\ra\Omega^1_{Y/k}$ is bijective.
Similarly,  $Y$ is a \emph{Calabi--Yau manifold} if it is simply-connected with $h^i(\O_Y)=0$ for $1<i<\dim(Y)$.

Under strong additional assumptions, analogues of the $T^1$-Lifting Theorem 
(\cite{Ekedahl; Shepherd-Barron 2003}, \cite{Schroeer 2003}) and the Decomposition Theorem \cite{Patakfalvi; Zdanowicz 2019} hold true.
In light of the liftability of abelian varieties (\cite{Katz 1981},  \cite{Norman; Oort 1980}) and K3-surfaces \cite{Deligne 1976},
it is natural to wonder
whether any such $Y$ admits a   lifting to characteristic zero.
This, however, turns out to be false already for Calabi--Yau threefolds.
The first example   was given by Hirokado \cite{Hirokado 1999} in characteristic $p=3$.
The second author \cite{Schroeer 2004} found further examples in characteristic $p=2,3$ using quotients
of the   \emph{Moret-Bailly pencil of supersingular abelian surfaces} \cite{Moret-Bailly 1981}.
Further examples in dimension three at certain bounded sets of primes where constructed by Schoen \cite{Schoen 2009},
Hirokado, Ito and Saito (\cite{Hirokado; Ito; Saito 2007} and \cite{Hirokado; Ito; Saito 2008}),
Cynk and van Straten \cite{Cynk; van Straten 2009}, and Cynk and Sch\"utt \cite{Cynk; Schuett 2012}.
Finally, Achinger and Zdanowicz (\cite{Achinger; Zdanowicz 2017}, \cite{Zdanowicz 2021}) produced for each prime $p\geq 5$ a  non-liftable Calabi--Yau manifold
of dimension $2p$, based on the failure of Kodaira Vanishing as observed by Totaro \cite{Totaro 2019}.

The goal of this paper is to generalize the Moret-Bailly pencil to higher dimensions, and establish
further non-liftability results:
Let $A=E_1\times\ldots\times E_g$ be a  product of supersingular elliptic curves.
Roughly speaking, the embeddings of the local unipotent group scheme $\alpha_p=\GG_a[F]$ into the abelian variety $A$
are parameterized by the projectivization $\PP^n=\PP(\liea)$ of the Lie algebra $\liea=\Lie(A)$, where $n+1=g$.
In fact, any inclusion $\O_\PS(-d)\subset\O_\PS^{\oplus(n+1)}$ that is locally a direct summand 
corresponds to a family $H\subset A\times\PS$ of such  finite local group schemes.
Setting $X=A\times\PS$, we then form the family of quotients $Y=X/H$, which 
comes with  a radical surjection $\epsilon:X\ra Y$ and an induced fibration  $\varphi:Y\ra\PS$.

We call the fibrations $\varphi:Y\ra\PS$  and also the smooth proper   schemes $Y$   
\emph{Moret-Bailly families}, since they are higher-dimensional analogs of the famous construction of a non-isotrivial family of abelian surfaces
over the projective line  \cite{Moret-Bailly 1981}; such families were already mentioned by Grothendieck in
\cite{Grothendieck 1966}, Remark 4.6.

The dimension of $Y$ is $2n+1=2g-1$. It is easy to compute the Betti numbers $b_i(Y)$, 
but the cohomological invariants $h^i(\O_Y)$ remain mysterious.
Using a result of Achinger \cite{Achinger 2015} on the splitting type of the Frobenius push-forward   on toric varieties,
one may express $h^i(\O_Y)$ via  lattice point counts.
 It turns out that the  canonical projection $\psi:Y\ra A^{(p)}=A/A[F]$ is the Albanese map,
and the dualizing sheaf $\omega_Y$ is the pullback of $\O_\PS(m)$, for the integer $m=d(p-1)-g$.
In particular, we have   $c_1=0$ if and only if $d(p-1)=g$, and in this case 
$\omega_Y=\O_Y$.  More generally, $\omega_Y$ is \emph{anti-nef}, which means that $(\omega_Y\cdot C)\leq 0$
for every integral curve $C\subset Y$,  if and only if $d(p-1)\leq g$.
The first main result of this paper is:

\begin{maintheorem}
(See Thm.\ \ref{no projective lift}.)
Suppose $d(p-1)\leq g$ and $p\geq 3$. Then the Moret-Bailly family  $Y$ does not projectively lift to characteristic zero.
\end{maintheorem}
 
\smallskip
In the boundary case  $d(p-1)=g$, this apparently gives the first examples of non-liftable manifolds  with $c_1=0$ 
that do not belong the class of abelian varieties, hyperk\"ahler manifolds, Calabi--Yau manifolds,
or products thereof.
To establish the result, we  assume that a   projective lifting $\foY\ra\Spec(R)$ exists, and use the existence  of relative
Hilbert schemes $\Hilb_{\foY/R}$ to show that properties of the Albanese map
$V\ra\Alb_{V/\CC}$ for the resulting complex fiber $V=\foY\otimes_R\CC$  
would contradict a recent result of  Cao  \cite{Cao 2019}. 
In the case where $d(p-1)=g$, the contradiction could also be derived from the Beauville--Bogomolov Decomposition Theorem 
for K\"ahler manifolds with $c_1=0$. \footnote{In a first version of this paper, we considered only the case $d(p-1)=g$. 
Ludvig Olsson  discovered that our result could be extended to the situation where the dual of $\omega_Y$ is nef, 
if one avails oneself of the result of   Cao. We are very grateful to him for sharing his insight with us and allowing us to make use of it in this paper.}
 
We strongly believe that Moret-Bailly families do not even formally lift to characteristic zero. So far we
are not able to show this, but we can prove the following,
which takes into account the \emph{sign involution} on $Y$, viewed as a family of abelian varieties parameterized by $\PS$:

\begin{maintheorem}
{\rm (See Thm. \ref{no formal lift})} 
Suppose that $d=1$, $p-1\leq g$ and $p\geq 3$. Then $Y$ together with its sign involution does not formally lift to characteristic zero.
\end{maintheorem}

We  show that any formal lifting $\foY$ to characteristic zero, in which the sign involution extends,   admits an ample sheaf.
Then Grothendieck's Existence Theorem gives a contradiction to Theorem A. This argument relies
on  Rim's  equivariant deformation theory \cite{Rim 1980} and its generalization to mixed characteristics \cite{Schroeer; Takayama 2018},
and a  computation of  \emph{weights} in the groups 
$H^2(Y,\O_Y)$ and $H^1(Y,\Theta_Y)$ for  the action of $G=\{\pm 1\}$
coming from the sign involution on $A$. 
In the first version of the article,  we overlooked some contribution to the weights,
which led to an    unjustified stronger assertion.
It turns out that extending the sign involution to infinitesimal deformations $\foY$ of $Y$
is actually the same as extending the morphism $\varphi:Y\ra\PS$ 
(see Proposition \ref{equivalences for lifting}). 

In the last section of this article, we consider the possibility of lifting 
Moret-Bailly families to the ring $W_2(k)$  of Witt vectors of length two.
The third main result is the following:

\begin{maintheorem}
{\rm (see Thm.\ \ref{no lift to witt})}
Fix $n\geq 2$ and $d\geq 1$. If $n\not\equiv 2$ modulo $4$, and the ground field $k$ is perfect,  then   
the Moret-Bailly family $Y$ does not lift to $W_2(k)$ for almost all primes $p>0$.
\end{maintheorem}

Note that in contrast to the previous results, here there are no assumption on the dualizing sheaf $\omega_Y$.
Of course there still may be deformations some finite flat extensions of $W_2(k)$.
In the special case $n=3$, we actually can show that $Y$ does not lift to $W_2(k)$ for all primes $p\geq 7$.

To prove Theorem C, we construct
certain ample sheaves  on the Moret-Bailly family $Y$ that violate  the conclusion 
of the Kodaira vanishing theorem.   By the results of Deligne  and Illusie \cite{Deligne; Illusie 1987}, 
this ensures   non-liftability to $W_2(k)$.
The computation relies on the Hirzebruch--Riemann--Roch Theorem, together with
an analysis of the Todd class $\td(\Theta_Y)$ and properties of the Bernoulli numbers $B_i\in\QQ$.

\medskip
The paper is organized as follows:
In Section \ref{Families group schemes} we collect some facts on families
of algebraic group schemes in characteristic $p>0$ and discuss the four-term complex
that describes certain infinitesimal quotients.
In Section \ref{Families abelian varieties} we apply this to families of abelian varieties.
The Moret-Bailly families $\varphi:Y\ra\PS$ are introduced in Section \ref{Moret-Bailly}, where
we compute the dualizing sheaf and Betti numbers.
Section \ref{Picard and albanese} contains a description of the Picard scheme and the  Albanese map.
Fibers of the Albanese map play a crucial role in Section
\ref{Non-existence projective}, where we prove that   Moret-Bailly families
with $c_1=0$ do  not projectively lift to characteristic zero.
Here the main ingredient are relative Hilbert schemes, and the Beauville--Bogomolov decomposition
over the complex numbers.
In Section \ref{Cohomology and points} we express the cohomology groups $H^i(Y,\O_Y)$ in terms of 
cohomology on $\PP^n$ for coefficient sheaves that involve Frobenius pullbacks and exterior powers.
This is used in Section \ref{Cohomology and weights} to compute weights in  
$H^2(Y,\O_Y)$ and $H^1(Y,\Theta_Y)$, which are the crucial obstruction groups for infinitesimal deformations.
Section \ref{Non-existence formal} contains the proof that   Moret-Bailly families $Y$ with $c_1=0$, together with the sign involution,
do not admit a formal lift to characteristic zero, by using equivariant deformation theory in mixed characteristics. 
In Section \ref{Non-existence witt}  we show that there are ample invertible sheaves on $Y$ 
which violate the conclusion of the Kodaira vanishing theorem when $n,d,p$ satisfy certain 
numerical conditions.

\begin{acknowledgement}
We wish to thank the referee for careful reading and valuable suggestions, in particular for pointing out a mistake  in the first version  of
the proof for Proposition \ref{relevant weights}.
The second author is grateful for the hospitality during his visit at Pembroke College Oxford,
where much of the work was carried out. 
The research was supported by the Deutsche Forschungsgemeinschaft via the grant  \emph{SCHR 671/6-1 Enriques-Mannigfaltigkeiten}. 
It was also conducted       in the framework of the   research training group
\emph{GRK 2240: Algebro-geometric Methods in Algebra, Arithmetic and Topology}.
We like to thank Laurent Moret-Bailly for useful remarks on the literature, and Emilian Zdanowicz
for pointing out the results in \cite{Achinger; Zdanowicz 2017} and \cite{Zdanowicz 2021}. 
Last but not least, we are grateful to Ludvig Olsson for his remarks and for his contribution to this paper (see the discussion above).
\end{acknowledgement}

%===========================================================
\section{Families of algebraic group schemes}
\mylabel{Families group schemes}

Let $S$ be a base scheme.
A \emph{family of  algebraic group schemes}  is a scheme   $X$, together with  a morphism $X\ra S$ that is flat and
of finite presentation, endowed with the structure of a relative group scheme. 
Here we could allow algebraic spaces as well. However, for the sake of exposition  we stay in the realm of schemes,
which is sufficient for our applications.

Let us   assume that the sheaf of K\"ahler differentials $\Omega^1_{X/S}$ is locally free.
Then the tangent sheaf $\Theta_{X/S}=\uHom(\Omega^1_{X/S},\O_X)$ is locally free as well.
The \emph{sheaf of Lie algebras}  $\Lie_{X/S}$ is the pullback of $\Theta_{X/S}$ along the  neutral
section $e:S\ra X$. This is a locally free sheaf, endowed with a Lie bracket,
such that the fibers $\idealg=\Lie_{X/S}\otimes\kappa(a)$, $a\in S$ become  Lie algebras over the residue fields $\kappa(a)$.

Now suppose that $S$ has characteristic $p>0$. Then the sheaf of Lie algebras
$\Lie_{X/S}$ acquires the \emph{$p$-map} as additional structure, such that the $\idealg=\Lie_{X/S}\otimes\kappa(a)$
become \emph{restricted Lie algebras} over  $\kappa(a)$.
Recall that the map
\begin{equation}
\label{gabriel demazure}
\Hom(X,Y)\lra \Hom(\Lie_{X/S},\Lie_{Y/S})
\end{equation}
is bijective provided that $X$ has  height at most one, according to 
\cite{SGA 3a}, Expos\'e $\text{VII}_A$, Theorem 7.2. 
Note that the hom set on the right comprises $\O_S$-linear maps compatible with Lie bracket and $p$-map, and that \emph{height at most one} means 
that the relative Frobenius  
$F:X\ra X^{(p)}$ is trivial. 
Here $X^\frob$ is the pullback of $X$ along the absolute Frobenius map 
$F_S:S\ra S$, and the morphism $F$ comes from the commutative diagram
$$
\begin{CD}
X	@>F_X>>	X\\
@VVV		@VVV\\
S	@>>F_S>	S
\end{CD}
$$
of absolute Frobenius maps. 
Moreover,  each 
sheaf $\shH$ of restricted Lie algebras  that is  locally free of finite rank arises from some 
family of algebraic group scheme $H$ of height at most one.
In fact, $H$ is the relative spectrum of the sheaf of algebras 
\begin{equation}
\label{sheaf of algebras}
\shA=\shHom(U^{[p]}(\shH),\O_S),
\end{equation}
where $U^{[p]}(\shH)$ is the quotient of the sheaf $U(\shH)$ of universal enveloping algebras
by some sheaf of ideals defined via the $p$-map, as explained in \cite{SGA 3a}, Expos\'e $\text{VII}_A$, Section 5,
compare also  \cite{Demazure; Gabriel 1970}, Chapter II, \S7, No.\ 4.

Let $\shH\subset\Lie_{X/S}$ be a subsheaf that is locally a direct summand,
and assume that $\shH$ is stable under  both Lie bracket and $p$-map. 
Let $H\ra S$ be the corresponding family of group schemes of height at most one,
with $\Lie_{H/S}=\shH$.  
We now consider the inclusion $H\subset X$ and form the resulting quotient $Y=X/H$.
Such a quotient exists as an algebraic space.  It is actually a  scheme, because the projection
$X\ra Y$ is a finite universal homeomorphism
(\cite{Olsson 2016}, Theorem 6.2.2). If $H$ is normal, $Y$ inherits the structure of a family of algebraic groups.

\begin{proposition}
\mylabel{structure morphism smooth}
The structure morphism $Y\ra S$ is   flat and of finite presentation.
Moreover, it is smooth provided that  $X\ra S$ is smooth.
\end{proposition}

\proof
The projection $\epsilon:X\ra Y$ is   faithfully flat and of finite presentation, 
because it is  a torsor with respect to $H\times Y$. The assertion now follows from fppf descent.
\qed

\medskip
In the special case $\shH=\Lie_{X/S}$ the group scheme $H$ coincides with the kernel $X[F]$ of the relative Frobenius 
map. In the general case, we thus have an $S$-morphism   $X/H\to X^\frob$. 

We now assume that $H=X[F]$, and furthermore that  $X$   is smooth. Then the same holds for $Y$, and the homomorphism $X\ra X^\frob$
is an epimorphism, such that  $X/H=X^\frob$.
In turn we obtain    an exact sequence
$0\ra \shH\ra \Lie_{X/S}\ra \Lie_{Y/S}$
of families of restricted Lie algebras. 
By assumption, the inclusion on the left is locally a direct summand, so 
the cokernel $\shK=\Lie_{X/S}/\shH$ is locally free. Since forming the quotient $Y=X/H$ commutes
with base-change, the inclusion $\shK\subset\Lie_{Y/S}$ is locally a direct summand. 
Let $K\subset Y$ be the corresponding family of group schemes of height at most one. 
The isomorphism theorem ensures  $Y/K=X/X[F]=X^\frob$.  
This gives  a commutative diagram:
$$
\begin{tikzcd}
0\ar[r]		&	\shH\ar[r]		&	\Lie_{X/S}\ar[r]\ar[dr]	&	\Lie_{Y/S}\ar[d]\ar[dr]\\
			&	0\ar[r]			&	\shH^\frob\ar[r]	&	\Lie_{X^\frob/S}\ar[r]	&	\Lie_{Y^\frob/S}
\end{tikzcd}
$$
The two diagonal maps vanish, hence the vertical map factors over $\shH^\frob$ and we obtain the \emph{four-term complex}
\begin{equation}
\label{lie algeba sequence}
0\lra \shH\lra \Lie_{X/S}\lra \Lie_{Y/S}\lra \shH^\frob \lra 0.
\end{equation}
Sequences like this already appear  in Ekedahl's work on foliations of smooth algebraic schemes
(\cite{Ekedahl 1987}, Corollary 3.4). 

\begin{theorem}
\mylabel{lie sequence exact}
The above complex of restricted Lie algebras is exact.
\end{theorem}

\proof
By construction, the complex is exact at all terms, with the possible exception of $\shH^\frob$.
Our task is thus to verify that $\Lie_{Y/S}\ra\shH^\frob$ is surjective.
By the Nakayama Lemma, it suffices to do so after tensoring with $\kappa(a)$, for $a\in S$.
Since the formation of quotients commutes with base-change, so does the formation of the complex.
It thus suffices
to treat the case that $S$ itself is the spectrum of a field.
Now the terms become finite-dimensional vector spaces.
The outer terms have the same dimension, and the same holds for the inner terms.
In turn, the rank of $\Lie_{Y/S}\ra\shH^\frob$ coincides with the dimension of $\shH^\frob$,
so the map in question must be surjective.
\qed

%===========================================================
\section{Families of abelian varieties}
\mylabel{Families abelian varieties}

We keep the assumption of the preceding section, and assume
now that $\psi:X\ra S$ is a family of   abelian varieties of relative dimension $g\geq 0$.
According to \cite{Artin 1969}, Theorem 7.3 the  relative Picard functor $\Pic{X/S}$ is representable by an algebraic space, and
$P=\Pic_{X/S}^0$  is   called the \emph{family of dual abelian varieties}.
The   sheaf of Lie algebras is given by $\Lie_{P/S}=R^1\psi_*(\O_X)$, with trivial bracket.
Note that by results of Raynaud, any family of abelian varieties has schematic total space, but is not necessarily projective
(\cite{Faltings; Chai 1990}, Theorem 1.9 and \cite{Raynaud 1970}, Chapter VII, 4.2, compare 
also the discussion in \cite{Laurent; Schroeer 2021}, Section 4).

Let $X\ra Y$ be a homomorphism between two families of abelian varieties. It must be flat, by  the fiber-wise criterion for flatness
(\cite{EGA IVc}, Theorem 11.3.10). Suppose also that $\dim(X_s)=\dim(Y_s)$ for all points $s\in S$. 
Then the kernel $H\subset X$ is a family of finite group schemes, and the \emph{Cartier dual} $\uHom(H,\GG_m)$ is another
family of finite group schemes. As explained in \cite{Oda 1969}, Corollary 1.3, this sits in an exact sequence
\begin{equation}
\label{first pic sequence}
0\lra\uHom(H,\GG_m)\lra \Pic^0_{Y/S}\lra \Pic^0_{X/S}\lra 0.
\end{equation}
Note that this can also be seen with the  identification $\Pic^0_{X/S}=\uExt^1(X,\GG_m)$ explained in \cite{SGA 7a}, Expos\'e VII, and
the long exact Ext sequence.

Now suppose that the base scheme $S$ has characteristic $p>0$.
Let $\shH\subset\Lie_{X/S}$ be a subsheaf that is locally a  direct summand,
and stable under the $p$-map, and $H\subset X$ the corresponding family of group schemes
of height at most one. Then the quotient $Y=X/H$ is again a family of $g$-dimensional abelian varieties,
and the projection $X\ra Y$ induces an exact sequence \eqref{first pic sequence}.
Likewise, $0\ra H'\ra Y\ra X^\frob\ra 0$ with   $H'=X[F]/H$ gives a short exact sequence
\begin{equation}
\label{second pic sequence}
0\lra\uHom(H',\GG_m)\lra \Pic^0_{X^\frob/S}\lra \Pic^0_{Y/S}\lra 0.
\end{equation}
Write $\varphi:Y\ra S$ and $\psi:X\ra S$ for the structure morphisms.

\begin{proposition}
\mylabel{four-term sequence}
Suppose the $p$-map   $\Lie_{X/S}\ra\Lie_{X/S}$ factors over the subsheaf $\shH$. 
Then  we have a four-term exact sequence
$$
0\lra \shK\lra R^1\psi_*(\O_X)^\frob\lra R^1\varphi_*(\O_Y)\lra \shK^{(p)}\lra 0,
$$
with the sheaf  $\shK=\shHom(\Lie_{X/S}/\shH,\O_S)$. Moreover, the sequence is natural with
respect to the inclusion $H\subset X$.
\end{proposition}

\proof
Consider the family $H'=X[F]/H$ of group schemes of height at most one.
Its  Cartier dual $K=\uHom(H',\GG_m)$ is a family of finite group schemes.
The latter have height at most one, by our  assumption on the $p$-map on $\Lie_{X/S}$.
We see from 
\eqref{sheaf of algebras} that $H'$ and $K$ are the Frobenius kernels 
for the respective vector schemes $\Spec(\Sym^\bullet (\shH'^\vee))$ and $\Spec(\Sym^\bullet( \shK^\vee))$,
where  $\shH'=\Lie_{H'/S}$ and $\shK=\Lie_{K/S}$.

By Lemma \ref{cartier and sheaf duality} below, the sheaf of Lie algebras $\shK$
coincides with the linear dual of  $\shH'=\Lie_{X/S}/\shH$. 
Our assertion in now  is a consequence of Theorem \ref{lie sequence exact}, applied to 
the short exact sequence \eqref{second pic sequence}.
The four-term sequence is natural with respect to $H\subset X$, because the same holds for
the two short exact sequences \eqref{first pic sequence} and \eqref{second pic sequence}.
\qed

\medskip
\newcommand{\LocLib}{{\text{\rm LocLib}}}
The preceding proof relies on the following  observations:
Let $(\LocLib/S)$ be the category of locally free sheaves of finite rank, and $(\Grp/S)$ be the category 
of families of algebraic group schemes.  Consider the functors
$$
\shE\longmapsto V\quadand \shE\longmapsto V^*,
$$
where $V=\Spec(\Sym^\bullet(\shE^\vee))$ and $V^*=\Spec(\Sym^\bullet(\shE))$ are the vector bundles
with $\Lie_{V/S}=\shE$ and $\Lie_{V^*/S}=\shE^\vee$, with trivial Lie brackets and $p$-maps.
Note that we follow Grothendieck's convention from \cite{EGA I}, Section 9.6.
The Frobenius kernels $G=V[F]$ and $G^*=V^*[F]$ are   families of finite local group schemes.

\begin{lemma}
\mylabel{cartier and sheaf duality}
The contravariant functors $\shE\mapsto G^*$ and $\shE\mapsto\uHom(G,\GG_m)$ are    naturally isomorphic.
In particular, there is an   identification 
$$
\Lie_{\uHom(G,\GG_m)/S}=\shHom (\Lie_{G/S},\O_S)
$$
of locally free sheaves that is natural in $G$.
\end{lemma}

\proof 
The natural identification  arises as follows:
Let $T=\Spec(A) $ be an affine  $S$-scheme,
and consider the resulting  $A$-module $E=\Gamma(T,\shE_T)$, which is finitely generated and projective.
According to \eqref{gabriel demazure}, the group of $A$-valued points in the Cartier dual $\uHom(G,\GG_m)$
is the set of linear maps 
$$
E=\Lie_{G/S}\otimes A\lra \Lie_{\GG_m/S}\otimes A=A
$$ 
that are
compatible with $p$-maps.  On the left-hand side, the $p$-map vanishes, 
whereas on the right hand side it is nothing but $\lambda\mapsto\lambda^p$.
So these linear maps can be seen as vectors in the dual $\Hom_A(E,A)$ annihilated by the relative Frobenius map.
The latter   coincide with the $A$-valued points of the Frobenius kernel $G^*=V^*[F]$.
\qed

\medskip
The following property of families of abelian varieties $\varphi:Y\ra S$ will be crucial for later 
computations:

\begin{proposition}
\mylabel{spectral sequence trivial}
For all $s\geq 0$, the higher direct image sheaves $R^s\varphi_*(\O_Y)$ are locally free, their formation commutes with base-change,
and the canonical maps $\Lambda^sR^1\varphi_*(\O_Y)\ra R^s\varphi_*(\O_Y)$ are bijective.
Moreover, the Leray--Serre spectral sequence
$$
E_2^{r,s}=H^r(S,R^s\varphi_*(\O_Y))\Longrightarrow H^{r+s}(Y,\O_Y)
$$
has trivial differentials on the $E_i$-page provided $p-1$ does not divide $ i-1$.
\end{proposition}

\proof
The first assertion is \cite{Berthelot; Breen; Messing 1982}, Proposition 2.5.2.
For the second assertion,   note that the differentials on the $i$-th page 
are certain additive maps 
$
d_i: E_i^{r,s}\ra E_i^{r+i,s-i+1}$,
which  are natural with respect to the family of abelian varieties $Y$. In particular, for each integer $n$ the
multiplication   on $Y$  induces an endomorphism $n^*$ on $\Lambda^sR^1\varphi_*(\O_Y)$. The latter is multiplication
by $n^s$. To check this it suffices to treat the case where    $S=\Spec(R)$ is affine and  $s=1$.
Then the families of abelian varieties form an additive category,
and $Y\mapsto H^1(Y,\O_Y)$ is a contravariant functor into the additive category  of all  $R$-modules,
with the property  $H^1(\O_{Y_1\times Y_2})=H^1(\O_{Y_1})\oplus H^1(\O_{Y_2})$.
Now recall that any functor $F:\catC\ra \catC'$ between additive categories that respects products
also respects the $\ZZ$-module structure on hom sets (\cite{Kashiwara; Schapira 2006}, Proposition 8.2.15),
and thus $n^*=n$.

Now suppose that there is an element $a\in E_i^{r,s}$ whose image $b\in E_i^{r+i,s-i+1}$ is non-zero.
The latter  can be seen as a basis vector inside an $\FF_p$-vector space.
Recall that the multiplicative group $\FF_p^\times $ is cyclic of order $p-1$. Choose an integer $n$
that generates $\FF_p^\times$. Since $p-1\nmid i-1$ we have $n^{i-1}\not\equiv 1$ modulo $p$.
One gets 
$$
n^sb=n^sd_i(a)=d_i(n^sa) = n^{s-i+1}d_i(a) = n^{s-i+1}b
$$
from the naturality of the Leray--Serre spectral sequence applied to $n^*$.
Comparing coefficients gives $n^{1-i}\equiv 1$ modulo $p$, contradiction.
\qed

%===========================================================
\section{Moret-Bailly families}
\mylabel{Moret-Bailly}

Let $k$ be a ground field of characteristic $p>0$. 
Recall that an abelian variety $A$ of dimension $g\geq 1$ is called \emph{superspecial}
if the Lie algebra $\idealg=\Lie(A)$ has trivial $p$-map.
For $g=1$ this means that $A=E$ is a supersingular elliptic curve.
Moreover, the products $A=E_1\times\ldots\times E_g$ of supersingular elliptic curves
are superspecial. Note that if $k$ is algebraically closed, the converse holds
(\cite{Oort 1975}, Theorem 2). If moreover $g\geq 2$,
the isomorphism class of $A$ does not depend on the factors
(\cite{Shioda 1978}, Theorem 3.5). We need the following well-known existence result:

\begin{lemma}
\mylabel{superspecial exist}
In each dimension $g\geq 1$  there is a superspecial abelian variety $A$.
\end{lemma}

\proof
Using the three Weierstra\ss{} equations 
$$
y^2+xy=x^3+ \frac{36}{1728-j}x + \frac{1}{1728-j}\quadand y^2+y=x^3\quadand y^2=x^3+x,
$$
one sees that each invariant $j\in k$ is attained by some elliptic curve
(compare \cite{Tate 1972}, Example on page 36).
So it suffices to show
that there are supersingular $j$-values in the prime field $k=\FF_p$.
For $p=2$ this is $j=0$. 
Suppose now $p\geq 3$.
Recall that an elliptic curve in Legendre form $E:y^2=x(x-1)(x-\lambda)$ is supersingular if and only
if $\lambda$ is a root of the Hasse polynomial $P(T)=\sum_{i=0}^m\binom{m}{i}^2T^i$,
where $m=(p-1)/2$. One may view  the spectrum of $k[\lambda]$ as 
the coarse moduli space for elliptic curves $E$ with  level structure $(\ZZ/2\ZZ)^2_k\subset E$. The group 
$G=\GL_2(\FF_2)$ acts freely via the level structures, and the ring of invariants $k[\lambda]^G=k[j]$
is the coarse moduli space for the Deligne--Mumford stack $\shM_{1,1}$. 
According to \cite{Brillhart; Morton 2004}, Theorem 1 the Hasse polynomial has at least one root over $k$ if and only
if $p\not\equiv 1$ modulo $4$; then a supersingular $E$ over $\FF_p$ already exists with   level structure,
and can be put in Legendre form.
Suppose now $p\equiv 1$ modulo $4$, and consider the subgroup $H\subset G$ generated interchanging the
2-division points with coordinate $x=0,1$. This is given by the change of coordinates $x=-x'+1$,
and induces  $\lambda\mapsto 1-\lambda$ on the coarse moduli space.
By loc.\ cit.\  the Hasse polynomial
viewed as element in the ring of invariants $k[\lambda]^H=k[ \lambda-\lambda^2]$
acquires  a root, which gives the desired supersingular $j$-value.
\qed

\medskip
Fix some integers $n,d\geq 1$ and choose a superspecial abelian variety $A$ of dimension $g=n+1$.
In turn, every non-zero vector in $\lieg=\Lie(A)$ gives an inclusion $\alpha_p\subset A$.
Choose an identification $\idealg=k^{n+1}$.
Set $X=A\times\PP^n$, and view this as the constant  family of abelian varieties
over $\PP^n=\Proj k[T_0,\ldots,T_n]=\PP(\lieg^\vee).$
Choose some homogeneous polynomials $Q_0,\ldots,Q_n$ of degree $d\geq 1$
without common zero on the projective space,
and consider the resulting inclusion 
$$
\O_S(-d) \subset \O_\PS^{\oplus n+1}=\lieg\otimes_k \O_\PS=\Lie_{X/\PS}.
$$
Let $H\subset X$ be the family of height-one group schemes    with $\Lie_{H/S}=\O_S(-d)$,
and form the resulting quotient $Y=X/H$. This is smooth and proper, of dimension $\dim(Y)=2n+1$, and with $h^0(\O_Y)=1$.
The inclusion $H\subset A[p]\times\PP^n$ induces a finite morphism $Y\ra (A/A[F])\times \PP^n$, hence $Y$ is   projective.
Write 
$$
\varphi:Y=X/H=(A\times\PS)/H\lra\PS
$$
for the structure morphism, which is  a family of supersingular abelian varieties
of dimension $g=n+1$, and $\epsilon:X\ra Y$ for the quotient map. 
 
We call $Y$ a \emph{Moret-Bailly family},
because the above generalizes the pencils in  \cite{Moret-Bailly 1981}, Part 2 to arbitrary dimensions 
(where the case $n=1$ and $d=1$ is considered). 
Note that the construction $Y=Y_{A,q}$
depends on the superspecial abelian variety $A$ and the finite flat morphism $q:\PS\ra\PS$
defined by the homogeneous polynomials $Q_i$, but we usually neglect this in notation.

The cokernel $\shE_d$ for the inclusion $\O_\PS(-d)\subset\lieg\otimes\O_\PS$ 
is locally free and sits in the short exact sequence
\begin{equation}
\label{defining sequence}
0\lra \O_\PS(-d)\lra  \O_\PS^{\oplus n+1} \lra \shE_d\lra 0,
\end{equation}
thus has $\det(\shE_d)=\O_\PS(d)$.
Note that for $Q_i=T_i$ this becomes the \emph{Euler sequence}
(compare \cite{Okonek; Schneider; Spindler 1980}, page 6), and   $\shE_1=\Theta_{\PS/k}(-1)$.
In general, we   have
$\shE_d=q^*(\Theta_{\PS/k}(-1))$, where $q:\PS\ra\PS$ is the morphism of degree $\deg(q)=d^n$
defined by the homogeneous polynomials $Q_i$.

Note that the Frobenius pullback  $\shE_{pd}=\shE_d^\frob$ is obtained by taking the $p$-powers  $Q_0^p,\ldots,Q_n^p$.
To simplify notation, we write  $\shE_{-d}=\uHom(\shE_d,\O_\PS)$ for the dual sheaves, and also
set $\shE_{-pd}=\shE_{-d}^\frob$. 
% Dualizing \eqref{} we obtain an exact sequence
% \begin{equation}
% \label{dual defining sequence}
% 0\lra \shE_{-d}\lra \O_\PS^{\oplus n+1} \lra \O_\PS(d)\lra 0.
% \end{equation}
% In fact, any inclusion of $\shE_{-d}$ into $\O_\PS^{\oplus n+1}$ that is locally a direct summand yields
% such a short exact sequence, which can be seen by taking determinants.

\begin{proposition}
\mylabel{lie and image}
The sheaf of Lie algebras and the first direct image are given by 
$$
\Lie_{Y/\PS}=\shE_d\oplus\O_\PS(-dp)\quadand 
R^1\varphi_*(\O_Y)=\O_\PS(d)\oplus\shE_{-pd}.
$$
Moreover, the $p$-map in $\Lie_{Y/\PS}$ 
is trivial on the first summand, and sends the second summand to the first, for all such splittings.
\end{proposition}

\proof
By Theorem \ref{lie sequence exact}, the sheaf of Lie algebras is an extension of  
$\O_\PS(-pd)$ by $\shE_d$.
All such extensions split: We have $\Ext^1(\O_\PS(-pd),\shE_d)=H^1(\PS,\shE_d(pd))$, and 
tensoring \eqref{defining sequence} with $\O_\PS(pd)$ yields an exact sequence
$$
H^1(\PS,\lieg\otimes_k\O_\PS(pd))\lra H^1(\PS,\shE_d(pd))\lra H^2(\PS,\O_\PS(pd-d)).
$$
The outer terms vanish, because the invertible sheaves $\O_\PS(pd)$ and $\O_\PS(pd-d)$ are ample.
Thus $\Lie_{Y/\PS}$ splits.

Next we verify  the assertion on the $p$-map. By construction, $\shE_d$ is a quotient
for $\Lie_{X/S}$, whereas $\O_\PS(-pd)$ is a subsheaf of $\Lie_{X^\frob/S}$.
Since $X=A\times\PS$ comes from a superspecial abelian variety $A$, the $p$-maps vanish
on these sheaves of Lie algebras, and the assertion follows.

It remains to analyze the first direct image. The coherent  sheaf 
$$
\shK=\uHom(\O_X^{\oplus n+1}/\O_\PS(-d),\O_\PS)=\uHom(\shE_d,\O_\PS)
$$
is nothing but $\shE_{-d}$, and its Frobenius pullback becomes $\shK^\frob=\shE_{-pd}$.
According to Proposition \ref{four-term sequence},
we have a four term exact sequence
$$
0\lra \shE_{-d}\lra \O_\PS^{\oplus n+1} \lra R^1\varphi_*(\O_Y)\lra \shE_{-pd}\lra0.
$$
The inclusion in the left is locally a direct summand, so its cokernel is isomorphic to $\O_\PS(d)$,
which follows by taking determinants.
Thus the direct image $R^1\varphi_*(\O_Y)$ is an extension of $\shE_{-pd}$ by $\O_\PS(d)$.
As in the preceding paragraph, one verifies that all such extensions split.
\qed

\medskip
This computation has the following consequence:

\begin{corollary}
\mylabel{field of definition}
For each rational point $a\in\PS$, there are only finitely many other rational points
$b\in\PP^n$ such that $\varphi^{-1}(a)\simeq\varphi^{-1}(b)$.
Moreover, any  field of definition $F$ for the generic fiber $\varphi^{-1}(\eta)$
has $\trdeg(F)=n$.
\end{corollary}

\proof
We may assume that $k$ is algebraically closed. Choose some odd prime $l\neq p$,
and some symplectic level structure $(\ZZ/l\ZZ)^{2g}\ra A$. This descents to a family of symplectic level structures for $Y$.
Let $\shA_{g,l}$ be the Artin stack of $g$-dimensional abelian varieties endowed with such a   structure.
This is actually an algebraic space (\cite{Faltings; Chai 1990}, Chapter IV, Corollary 2 on page 131) that is separated and of finite type.
Our Moret-Bailly family corresponds to a morphism $h:\PS\ra\shA_{g,l}$. It suffices  to check that $h$
is  quasi-finite, because every abelian variety has only finitely many such level structures.
Suppose it  is not quasi-finite. Then there is an integral curve $C\subset\PS$ that maps to a closed point.
Then the restriction $Y_C=Y\times_\PS C$ is isomorphic to $B_C=B\times C$ for some abelian variety $B$.
Passing to the sheaf of Lie algebras we get  $(\shE_d \oplus\O_\PS(-dp))_C\simeq\O_C^{\oplus n+1}$ as coherent sheaves,
in contradiction to the Krull--Schmidt Theorem for coherent sheaves (\cite{Atiyah 1956}, Theorem 1).
The assertion on the field of definition is proven in an analogous way.
\qed

\medskip
To simplify notation, we now  write $\O_Y(m)$ for the preimages under the structure morphism $\varphi:Y\ra \PS$
of the invertible sheaves $\O_{\PP^n}(m)$.

\begin{corollary}
\mylabel{dualizing sheaf}
The   dualizing sheaf takes the form  $\omega_Y=\O_Y(m)$ for the integer $m=d(p-1)-(n+1)$.
In particular, $c_1=0$ holds if and only if $d(p-1)=n+1$, and in this case we  actually have
$\omega_Y=\O_Y$.
\end{corollary}

\proof
Since $Y$ is smooth, the dualizing sheaf is $\omega_Y=\det(\Omega^1_{Y/k})$.
We have a short exact sequence $0\ra \varphi^*(\Omega^1_{\PS/k})\ra \Omega^1_{Y/k}\ra\Omega^1_{Y/\PS}\ra 0$.
The cokernel is isomorphic to the preimage of the dual for $\Lie_{Y/\PS}$, because $\varphi:Y\ra\PS$
is a family of smooth algebraic group schemes.
The sheaf of Lie algebras equals $\shE_d\oplus \O_\PS(-pd)$, and we have $\det(\shE_d)=\O_\PS(d)$.
Furthermore $\det(\Omega^1_\PS)=\O_\PS(-n-1)$.
Combining all this we obtain  $\omega_Y=\O_Y(m)$ for the integer $m=-n-1 -d + pd$.
\qed

\medskip
It follows that the Kodaira dimension takes the values
$\kod(Y) \in\{-\infty,0,n\}$
depending on the sign of the integer $m=d(p-1)-(n+1)$.
Moreover, in  the canonical model $Z=\Proj\bigoplus_{t\geq 0} H^0(Y,\omega^{\otimes t})$ 
is given by the respective schemes $Z=\varnothing, \PP^0, \PP^n$. Analogous statements
hold  for the anticanonical models.

It is   easy to determine the Betti numbers $b_i\geq 0$,  defined as the ranks
of the $l$-adic cohomology groups $H^i(\bar{Y},\ZZ_l(i))=\invlim_{\nu\geq 0} H^i(\bar{Y},\mu_{l^\nu}^{\otimes i})$,
where  $\bar{Y}=Y\otimes k^\alg$ is the  base-change to some algebraic closure, and $l>0$ is a prime
different from $p$.

\begin{proposition}
\mylabel{l-adic cohomology}
The $l$-adic cohomology groups $H^i(\bar{Y},\ZZ_l(i))$ are free of rank  
$b_i= \sum_j  \binom{2n+2}{i-j}$,
where the sum runs over all even   $j\geq0$. 
In particular, we have $b_1=2n+2$ and $b_2=2n^2+3n+2$ and $b_{2n+1}=2^{2n+1}$.
\end{proposition}

\proof
We may assume that $k$ is algebraically closed.
The quotient map  $\epsilon:X\ra Y$ is a finite universal homeomorphism, so the $l$-adic
cohomology groups for $Y=X/H$ and $X=A\times\PS$ coincide.
Taking cohomology with coefficients in $R=\ZZ/l^\nu\ZZ$, we have
$H^\bullet(\PP^n) =R[h]/(h^{n+1})$  and  $H^\bullet(A) = \Lambda^\bullet H^1(A)$.
These $R$-modules are free, where the generator $h$ has degree two,  and $H^1(A)$ is of rank $2n+2$.
In turn, 
$$
H^i(X) = \bigoplus_{j} H^j(\PP^n)\otimes_R H^{i-j}(A)
$$
by the K\"unneth Formula (\cite{SGA 4c}, Expos\'e XVII, Theorem 5.4.3) and the assertion 
on $H^i(Y,\ZZ_l(i))$ is a direct consequence. The values $b_1$ and $b_2$ follow immediately.
In middle degree, we get $b_{2n+1}=\sum_s  \binom{2n+2}{s}$, where the sum runs over all odd $s$.
This sum is  half of  $(1+1)^{2n+2} - (1-1)^{2n+2} = 2^{2n+2}$.
\qed

\medskip
Note that since $b_{2n+1}\neq 0$,  the method introduced by Hirokado \cite{Hirokado 1999} to establish non-liftability  
apparently does not apply.

%===========================================================
\section{The Picard scheme and the Albanese map}
\mylabel{Picard and albanese}

Keep the notation from the previous section, such that $Y=(A\times\PS)/H$
is a Moret-Bailly family formed with some superspecial abelian variety $A$
of dimension $g=n+1$ 
and some homogeneous polynomials $Q_0,\ldots,Q_n$ of degree $d\geq 1$, without common zero on $\PP^n$.
Let $\varphi:Y\ra\PS$ be the structure morphism.
We now examine  the Picard scheme for $Y$.  Recall that its Lie algebra is the cohomology group $H^1(Y,\O_Y)$.

\begin{proposition}
\mylabel{dimension pic}
The Picard scheme $\Pic_{Y/k}$ has dimension $n+1$, and furthermore  $h^1(\O_Y)=\binom{n+d}{d}$.
In particular, the Picard scheme is smooth if and only if $d=1$.
\end{proposition}

\proof
We may assume that $k$ is algebraically closed. To compute the dimension $d\geq 0$ of the Picard scheme, choose a prime $l\neq p$
that does not divide the order of the torsion part in $\NS(Y)$.
The Kummer sequence $0\ra\mu_l\ra\GG_m\ra\GG_m\ra 0$ implies that $2d=b_1$.
According to Proposition \ref{l-adic cohomology}, we have $b_1=2n+2$.

The Leray--Serre spectral sequence for $\varphi:Y\ra\PP^n$ gives an exact sequence
$$
0\lra H^1(\PS,\O_\PS)\ra H^1(Y,\O_Y)\lra H^0(\PS, R^1\varphi_*(\O_Y)) \lra H^2(\PS,\O_\PS).
$$
The outer terms vanish, and we merely have to compute the global sections of $R^1\varphi_*(\O_Y)=\O_\PS(d)\oplus\shE_{-pd}$.
The first summand contributes $h^0(\O_\PS(d))=\binom{n+d}{d}$.
It remains to check that the sheaf $\shE_{-pd}$ has no non-zero global sections.
Dualizing the short exact sequence $0\ra\O_\PS(-pd)\ra\O_\PS^{\oplus n+1}\ra\shE_{pd}\ra 0$, we get 
an exact sequence
$$
0\lra H^0(\PS,\shE_{-pd})\lra H^0(\PS,\O_\PS^{\oplus n+1})\lra H^0(\PS,\O_\PS(pd)).
$$
The map on the right is given by the homogeneous polynomials $Q_0^p,\ldots,Q_n^p$.
These are linearly independent,  so the map must be injective.
It follows that $H^0(\PS,\shE_{-pd})=0$. This shows $h^1(\O_Y)=\binom{n+d}{d}$.
\qed

\medskip
Set $X=A\times\PS$.
The family $H\subset X$ of finite group schemes of height one sits inside the constant family
$X[F]=A[F]\times\PS$.
In turn,  we get an induced homomorphism $Y\ra A^\frob\times\PS$ between families of abelian varieties,
and write $\psi:Y\ra A^\frob$
for the composition with the projection.
 
\begin{proposition}
\mylabel{properties albanese}
The morphism $\psi:Y\ra A^\frob$ is flat, and every geometric fiber is non-reduced, with reduction isomorphic
to the projective $n$-space. Moreover, the canonical map $\O_{A^\frob}\ra\psi_*(\O_Y)$ is bijective.
\end{proposition}

\proof
The quotient map $\epsilon:X\ra Y$ if faithfully flat, and so is the composition $\psi\circ\epsilon:X=A\times\PS\ra A^\frob$.
By descent,  $\psi:Y\ra A^\frob$ must be flat. 

The fiber $Z=\psi^{-1}(0)$ over the origin is the family $(A[F]\times\PS)/H$ of height-one group schemes,
with $\Lie_{Z/\PS}=\shE_d$. We claim that $h^0(\O_Z)=1$.
The universal enveloping algebra is $U(\shE_d)=\Sym^\bullet(\shE_d)$,
and the restricted quotient becomes $U^{[p]}(\shE_d)=\Sym^\bullet(\shE_d)/\shE_d^\frob\Sym^\bullet(\shE_d)$.
According to \eqref{sheaf of algebras},  $Z$ is the relative spectrum of the corresponding sheaf of Hopf algebras
$$
\shA=\shHom( U^{[p]}(\shE_d),\O_{\PS})\subset \shHom(\Sym^\bullet(\shE_d),\O_\PS).
$$
The term on the right is the   product of the coherent sheaves
$\shHom(\Sym^i(\shE_d),\O_\PS)$, and the term on the left is already contained in the corresponding
sum. The summands are  \emph{divided powers}
$\Gamma^i(\shF)=\Sym^i(\shF^\vee)^\vee$, for the dual sheaf $\shF=\shE_d^\vee$.
It thus suffices to verify that the   divided powers have no
non-zero global sections for $i\geq 1$. We proceed by induction.
The case $i=1$ was already treated in the proof for Proposition \ref{dimension pic}.
Now suppose $i\geq 2$, and that the assertion is true for $i-1$.
The surjection $\O_\PS^{\oplus n+1}\ra\shE_d$ induces a canonical surjection
$\Sym^{i-1}(\shE_d)\otimes\O_\PS^{\oplus n+1}\ra\Sym^i(\shE_d)$.
Dualizing the latter gives an inclusion $\Gamma^i(\shF)\subset \Gamma^{i-1}(\shF)\otimes\O_\PS^{\oplus n+1}$,
which completes the induction. Note that the inclusion 
is a  piece from the Eagon--Northcott complex \cite{Eagon; Northcott 1962}.
 Summing up, this establishes $h^0(\O_Z)=1$.

Now consider the fiber over a geometric point $\bar{b}:\Spec(\Omega)\ra A^\frob$, with image $b\in A^\frob$.
Making a base-change, it suffices to treat the case that $k=\Omega$ is algebraically closed and
that $b\in A^\frob$ is rational. Then $b=a^\frob$ for some rational point $a\in A$, and
  $\psi^{-1}(b) = ((a+A[F])\times \PS)/H$. This is isomorphic to $\psi^{-1}(0)=(A[F]\times\PS)/H$ via translation by $a$,
so the fiber is non-reduced, with reduction isomorphic to $\PS$.

It remains to compute $\psi_*(\O_Y)$. We just saw that the function $b\mapsto h^0(\O_{Y_b})=1$ is  constant on the reduced
scheme $A^\frob$. It follows that the direct image sheaf is locally free of rank one,
hence the canonical map $\O_{A^\frob}\ra\psi_*(\O_Y)$ is bijective.
\qed
 
\medskip
As explained by Serre \cite{Serre 1960}, there is a morphism $Y\ra V$ to some abelian variety $V$
such that every other such morphism $Y\ra V'$ arises via composition with some unique   $V\ra V'$.
Note that the latter usually does not respect the origin.
This  $V=\Alb_{Y/k}$ is called the \emph{Albanese variety}, and $Y\ra\Alb_{Y/k}$ is the \emph{Albanese map}.
See \cite{Laurent; Schroeer 2021} for a relative theory.

\begin{proposition}
\mylabel{albanese map}
The morphism $\psi:Y\ra A^{(p)}$ is the Albanese map.  
\end{proposition}

\proof
Let $f:Y\ra B$ be  a morphism into another abelian variety. We have to show that this map factors uniquely
over $\psi:Y\ra A^\frob$.
The composition $f\circ \epsilon:X\ra B$ factors over the projection $\pr_1:X=A\times\PP^n\ra A$. This gives a commutative
diagram
$$
\begin{CD}
X	@>\pr_1>> 	A\\
@V\epsilon VV		@VVh V\\
Y 	@>>f> 	B.
\end{CD}
$$
Replacing $f$ by the composition of a translation of $B$ with $f$, we may assume that $h:A\ra B$ respects origins, hence is a homomorphism.
For each rational point $a\in\PP^n$, the fiber $H_a\subset X_a=A$ is a copy of $\alpha_p$
whose schematic image in $Y$ and hence also in $B$ is a rational point $b\in B$.
The inclusion $\O_\PS(-d)\subset\O_\PS^{\oplus n+1}$ defining the family of subgroup schemes $H\subset X$ 
is locally a direct summand, hence the canonical map  $\bigcup_a H_a \ra A$ from the disjoint union
has schematic image $A[F]$. In the commutative diagram
$$
\begin{tikzcd} 
\bigcup_a H_a\ar[r,"\can"]\ar[dr]   & H \ar[r,"\pr_1"]\ar[d,"0"']        & A\ar[d,"h"]\\
                                    & Y \ar[r,"f"']        & B
\end{tikzcd}
$$
the composition $h\circ\pr_1\circ\operatorname{can}$ factors over the origin $0\in B$,
and the schematic image of $\pr_1\circ\operatorname{can}$ is $A[F]$.
Thus $h:A\ra B$ factors over $A/A[F]=A^\frob$. 
The factorization is unique, because the composition $X\ra A\ra A/A[F] $ is faithfully flat, hence
an epimorphism.
\qed

%===========================================================
\section{Non-existence of projective liftings }
\mylabel{Non-existence projective}

We keep the assumptions of the preceding section, and furthermore assume that $n\geq 2$ and $d(p-1)\leq n+1$. These assumptions on $n,d,p$ are made only in this section. 
In this situation our Moret-Bailly family $Y=(A\times\PS)/H$
has dimension $\dim(Y)=2n+1\geq 5$ and $\omega_Y$ is anti-nef. Recall that one says  that $Y$ \emph{projectively lifts to characteristic
zero} if there is a   local noetherian ring $R$ with residue field $k=R/\maxid_R$ such that the canonical map $\ZZ\ra R$ is injective,
together with a projective flat morphism $\foY\ra\Spec(R)$ whose closed fiber is isomorphic to $Y$.
Note that one may assume that $R$ is also complete and one-dimensional.
The goal of this section is to establish the following:

\begin{theorem}
\mylabel{no projective lift}
The scheme $Y$ does not projectively lift  to characteristic zero.
\end{theorem}

The proof requires some preparation and is given at the end of this section.
It suffices to treat the case that $k$ is algebraically closed.
Seeking a contradiction, we assume that $Y$ projectively lifts.
Then there is a complete discrete valuation ring $R$
with residue field $R/\maxid_R=k$ whose field of fractions $F=\Frac(R)$ has characteristic zero,
together with   a proper flat morphism $\upsilon:\foY\ra\Spec(R)$
with closed fiber $Y=\foY\otimes_Rk$. 
Write $V=\foY\otimes_R F$ for the generic fiber, which is a smooth proper scheme with $h^0(\O_V)=1$.
 
We start by examining the Picard scheme of $V$.
The component of the origin $P=\Pic^0_{V/F}$ is an abelian variety.
Consider the dual abelian variety $\Pic^0_{P/F}$.
After passing to a finite extension of $R$, we may assume that the structure morphism
$\foY\ra\Spec(R)$ admits a section. In particular, the generic fiber $V$ contains a rational
point. Then there is a Poincar\'e sheaf $\shP$ on $V\times P$, and we may assume that it is numerically
trivial on the fibers of the first projection. As explained in \cite{FGA VI}, Theorem 3.3, the resulting 
$$
\Psi:V\lra \Pic^0_{P/F},\quad v\longmapsto [\shP|\{v\}\times P]
$$ 
is the Albanese map, and we write $\Alb_{V/F}=\Pic^0_{P/F}$.

According to Proposition \ref{albanese map}, the   composition  of the quotient map with the projection 
$A\times\PP^n\ra A^\frob$ induces the Albanese map $\psi:Y\ra A^\frob$, and the reduced preimage of the origin  
$$
Z=\psi^{-1}(0)_\red=\PP^n_k.
$$
is a copy of the projective $n$-space. 
We now  exploit the existence of the relative Hilbert scheme $\Hilb_{\foY/R}$, which parameterizes flat families
of closed subschemes \cite{FGA IV}, and 
regard the closed subscheme $Z\subset Y$ as    a $k$-valued point $\xi =[Z]$ in the relative Hilbert scheme.

\begin{proposition}
\mylabel{hilbert smooth}
The structure morphism $\Hilb_{\foY/R}\ra \Spec(R)$ is smooth near $\xi $.
\end{proposition}

\proof
Let $\shI\subset\O_Y$ be the sheaf of  ideal for the closed subscheme $Z \subset Y$.
According to \cite{FGA IV}, Corollary 5.4 it suffices to check that the obstruction group $\Ext^1(\shI/\shI^2,\O_Z)$
vanishes. Since the inclusion $Z\subset Y$ is a regular embedding,
the conormal sheaf $\shN=\shI/\shI^2$ is locally free, and the obstruction group becomes $H^1(Z,\shN^\vee)$.

To compute the sheaf $\shN$, we consider the commutative diagram
\begin{equation}
\label{hilbert diagram}
\begin{tikzcd} 
Z\ar[r]\ar[dr,"\simeq"']	&	Y\ar[d,"\varphi"]\ar[r]		& \Spec(k)\\\
						&	\PS\ar[d]\\
						&	\Spec(k)
\end{tikzcd}
\end{equation}
where the arrows are either smooth or regular embeddings. The vertical part yields a short exact sequence
$0\ra \varphi^*\Omega^1_{\PS/k}\ra \Omega^1_{Y/k}\ra \Omega^1_{Y/\PS}\ra 0$.
The  kernel is locally a direct summand, so the restriction 
$$
0\lra \varphi^*\Omega^1_{\PS/k}|Z\lra \Omega^1_{Y/k}|Z\lra \Omega^1_{Y/\PS}|Z\lra 0,
$$
remains exact. The horizontal part yields another short exact sequence
$$
0\lra \shI/\shI^2\lra \Omega^1_{Y/k}|Z\lra \Omega^1_{Z/k}\lra 0.
$$
Using the commutativity of \eqref{hilbert diagram}, we see that the inclusion of  $\varphi^*\Omega^1_{\PS/k}|Z$
splits the above extension, so the projection $\shI/\shI^2\ra \Omega^1_{Y/\PS}|Z$ is bijective.
With the identification $Z=\PP^n$ and Proposition \ref{lie and image}  we conclude $\shN^\vee=\shE_d\oplus\O_\PS(-pd)$.

One now easily checks that $H^1(Z,\shN^\vee)=0$.
Indeed, the short exact sequence $0\ra\O_\PS(-d)\ra\O_\PS^{\oplus n+1}\ra\shE_d\ra 0$
yields an exact sequence
$$
H^1(\PS,\O_\PS^{\oplus n+1})\lra H^1(\PS,\shE_d)\ra H^2(\PP^n,\O_\PS(-d))
$$
The term on the left vanishes. This also holds for the term on the right
for $n\neq 2$. Since $n=d(p-1)\geq d$ this also remains true for $n=2$.
Furthermore, we have $H^1(\PS,\O_\PS(-pd))=0$, because $n\geq 2$.
\qed

\medskip
By Hensel's Lemma, the relative Hilbert scheme admits an $R$-valued point passing
through $\xi\in\Hilb_{\foY/R}$. Let $\foZ\subset\foY$ be the corresponding flat family
of closed subschemes, with closed fiber $\foZ\otimes_Rk=Z=\PP^n_k$.

\begin{proposition}
\mylabel{lifting projective space}
The generic fiber $\foZ\otimes_RF\subset\foY\otimes_RF=V$ is isomorphic to $\PP^n_F$, and  
must be contained in some fiber of the Albanese map $\Psi:V\ra\Alb_{V/F}$.
\end{proposition}
 
\proof
Using  $H^1(\PS,\Theta_\PS)=0$ we inductively construct compatible
isomorphisms $\PP^n\otimes{R/\maxid_R^{n+1}}\ra \foZ\otimes_RR/\maxid_R^{n+1}$.
Grothendieck's Existence Theorem gives $\PP^n_R\simeq\foZ$. Since every morphism
from the projective line to an abelian variety is constant, 
the scheme $\foZ_F$ must be contained in some fiber of $\Psi:V\ra\Alb_{V/F}$. 
\qed

\medskip
\emph{Proof of Theorem \ref{no projective lift}:} Note first that by applying the  Specialization Theorem (\cite{SGA 4c}, Expos\'e XVI, Corollary 2.2) to the smooth proper morphism
$\foY\ra\Spec(R)$, we may conclude that the Betti numbers of the closed and generic fiber coincide.
In combination with Proposition \ref{l-adic cohomology} this shows that $b_1(V)=2(n+1)$.
In turn, the abelian variety $\Pic^0_{V/F}$ has dimension   $n+1$, and the same holds
for the Albanese variety  $\Alb_{V/F}$. Applying \cite{Moriwaki 1992}, Corollary 8 to the dual of the relative dualizing sheaf of $\foY$, we also see that the dualizing sheaf $\omega_V$ of $V$ is nef. The main result of \cite{Cao 2019} now 
implies that the Albanese map $\Psi:V\ra\Alb_{V/F}$ is smooth. In particular, the fibres of 
$\Psi$ are equidimensional of dimension $n=2n+1-(n+1).$
By Proposition \ref{lifting projective space}, there is a closed point $\lambda\in\Alb_{V/F}$ such that the fiber $V_\lambda=\Psi^{-1}(\lambda)$ of the Albanese map $\Psi$ 
contains the generic fiber   $\foZ\otimes_R F$ of the family of subschemes $\foZ\subset\foY$. Furthermore, by the same proposition, we have $\foZ\otimes_R F\simeq\PP^n_F$. 
The  inclusion $\foZ\otimes_R F\subset V_\lambda$ must be a connected component, because $\dim(\PP^n_F)=\dim(V_\lambda)=n$
and $V_\lambda$ is smooth. Finally, the 
conormal bundle of $\foZ\otimes_R F$ in $V$ is trivial, because $\Psi$ is flat (since it is smooth). Using Corollary \ref{dualizing sheaf}, we thus obtain an isomorphism $$\omega_{\PP^n_F}\simeq\omega_V|_{\PP^n_F}\simeq 
\O(d(p-1)-n-1).$$ This is a contradiction, since $\omega_{\PP^n_F}\simeq\O(-n-1)$. 
\qed

%===========================================================
\section{Cohomology and lattice points}
\mylabel{Cohomology and points}

\newcommand{\Syz}{\operatorname{Syz}}
\newcommand{\define}{{\rm\text{def}}}
We now assume that the polynomials $Q_0,\ldots,Q_n\in k[T_0,\ldots,T_n]$
used to define  our Moret-Bailly family $Y=(A\times\PP^n)/H$
have degree $d=1$. After applying an automorphism of the projective $n$-space, we reduce to  the situation $Q_i=T_i$.
Tensoring the Euler sequence
$$
0\lra \O_\PS(-1)\lra \O_\PS^{\oplus n+1}\lra \shE_1 \lra 1
$$
with $\O_\PS(1)$ shows $\shE_1 =\Theta_\PS(-1)$, and the dual becomes
$\shE_{-1}=\Omega^1_\PS(1)$. We now use the Bott Formula 
$$
h^s(\Omega^r_\PS(l))= \begin{cases}
\binom{n+l-r}{n} \binom{l-1}{r}			& \text{if $s=0$ and $0\leq r\leq n$;}\\
\binom{r-l}{r}\binom{-l-1}{n-r}			& \text{if $s=n$ and $0\leq r\leq n$;}\\
1					& \text{if $0\leq r=s\leq n$ and $l=0$;}\\
0					& \text{else}
\end{cases}
$$
to understand  the groups $H^i(Y,\O_Y)$ better (compare \cite{Bott 1957}, Section 4.
For an algebraic proof, see \cite{Dolgachev 1981}, Section 2.3.2 or \cite{Huang 2001}, Section 4).
Note that the binomial coefficient  $\binom{x}{n}=x(x-1)\ldots(x-r+1)/n!$ 
is defined for natural numbers $n$ and     ring elements $x\in R$ whenever the denominator $n!$ is invertible.

Our computation crucially relies on the splitting type of the Frobenius push-forward
$F_*(\O_\PS)$, which indeed splits by the Horrocks Criterion 
(see   \cite{Okonek; Schneider; Spindler 1980}, Section 2.3).
Understanding the splitting type involves   seemingly innocent lattice point counts: 
For each integer $t$ and $l$, define multiplicities
$\mu_{t,l}\geq 0$ as the number of lattice points $(l_0,\ldots,l_n)\in\ZZ^{n+1}$
contained in the   polytope $P_{t,l}\subset \RR^{n+1}$ defined by 
\begin{equation}
\label{polytope}
l_0+\ldots+l_n =t-pl 	\quadand 	0\leq l_0,\ldots,l_n\leq p-1.
\end{equation}
This is  the intersection of an affine hyperplane with a hypercube.
Clearly, the polytope is non-empty if and only if $0\leq t-pl\leq (n+1)(p-1)$, and we have the
recursion formula
\begin{equation}
\label{recursion}
\mu_{t+p, l+1}=\mu_{t,l}.
\end{equation}
Note that the  $\mu_{t,l}\geq 0$ also  depends on $n$ and $p$, but we neglect this dependence in notation.
Applying a result of Achinger (\cite{Achinger 2015},  Theorem 2.1) to the toric variety $\PP^n$, we get:

\begin{proposition}
\mylabel{frobenius push-forward}
The Frobenius push-forward $F_*(\O_\PS(t))$ splits as a sum of invertible sheaves, and the summand  
$\O_\PS(l)$ appears with multiplicity $\mu_{t,l}\geq 0$.
In other words,    the splitting type is of the form
$$
(\underbrace{a,\ldots,a}_{\mu_{t,a}},\cdots,\underbrace{b,\ldots,b}_{\mu_{t,b}})
$$
starting with $a= \lceil(t-(n+1)(p-1))/p \rceil$and ending with $b=\lfloor t/p\rfloor$ .
\end{proposition}

In preparation for our analysis of $H^i(Y,\O_Y)$,
we  now express the cohomological invariants  of certain locally free sheaves on $\PP^n$ in terms of   lattice points
and binomial coefficients:

\begin{proposition}
\mylabel{exterior and frobenius}
The cohomological invariants of the locally free sheaf $\shF_{r,t}=\Lambda^r(F^*(\Omega^1_\PS(1)))\otimes\O_\PS(t)$
are given by the formula
$$
h^s(\shF_{r,t})= 
\begin{cases}
\sum_l \mu_{t,l} \binom{n+l}{n} \binom{l+r-1}{r}		& \text{if $s=0$;}\\
\sum_l \mu_{t,l}\binom{-l}{r}\binom{-r-l-1}{n-r}		& \text{if $s=n$;}\\
\mu_{t,-r}											& \text{if $s=r$;}\\
0														& \text{else.}
\end{cases}
$$
In particular, we have  $h^0(\shF_{r,t})=0$ for $0\leq t\leq p-1$, and $h^s(\shF_{0,0})=0$ for $s\geq 1$.
\end{proposition}

\proof
 Set $\shF=\shF_{r,t}$.
We have $H^s(\PS,\shF) = H^s(\PS,F_*\shF)$ because the Frobenius map is affine, and  the projection formula gives
$F_*\shF = \Omega^r_\PS(r)\otimes F_*(\O_\PS(t))$.
Now combine the Bott formula for $h^s(\Omega^r_\PS(l))$ and Achinger's description of $F_*(\O_\PS(t))$
to get the general formula for $h^s(\shF)$.
In the special case $s=0$ this reduces to
$$
h^0(\shF_{r,t})=\sum_l \mu_{t,l}\binom{n+l}{n}\binom{r+l-1}{r}.
$$
The second binomial coefficient   vanishes for $l\leq 0$. For $l\geq 1$ and $t\leq p-1$
the multiplicity $\mu_{t,l}$ is zero, because the polytope $P_{t,l}\subset\RR^{n+1}$ in  \eqref{polytope} becomes empty.
Finally, we have  $\shF_{0,0}=\O_\PS$ and thus $h^s(\shF_{0,0})=0$ in all degrees $s\geq 1$.
\qed

\medskip
We now can express  the cohomology groups of our Moret-Bailly family as follows:

\begin{proposition}
\mylabel{hodge numbers}
Suppose that $p\geq n+1$. For every degree $i\geq 0$, the Leray--Serre spectral sequence
for $\varphi:Y\ra\PS$ gives a natural identification
\begin{equation}
\label{natural identification}
H^i(Y,\O_Y) = \begin{cases}
H^j(\PP^n, \Lambda^j (\shE_{-p}	)	) 				& \text{if $i=2j$ is even;}\\
H^j(\PP^n, \Lambda^j(\shE_{-p}) \otimes\O_\PS(1))  	& \text{if $i=2j+1$ is odd.}
\end{cases}
\end{equation}
Moreover, the dimension of these vector spaces are given by the formula
$$
h^i(\O_Y) =\begin{cases}
\mu_{0,-j} & \text{if $i=2j$ is even;}\\
\mu_{1,-j} & \text{if $i=2j+1$ is odd,}
\end{cases}
$$
where   $\mu_{t,l}\geq 0$ is  the number of lattice points in the polytope $P_{t,l}\subset\RR^{n+1}$ as in \eqref{polytope}.
\end{proposition}

\proof
The assertion indeed holds for $i=0$, because $\varphi_*(\O_Y)=\O_\PS=\Lambda^0(\shE_{-p})$,
and  the only lattice point in $P_{0,0}$ has coordinates  $l_0=\ldots = l_n=0$.
Suppose from now on that  $i\geq 1$. 

The Leray-Serre spectral sequence is
$E_2^{i,j}=H^i(\PS,R^j\varphi_*\O_Y)\Rightarrow H^{i+j}(Y,\O_Y)$.
For dimension reasons, the differentials $d_r:E_r^{i,j}\ra E_r^{i+r,j-r+1}$ vanish on
the $E_r$-pages whenever $r\geq n+1$. Since $p-1\geq n$, they must also vanish on the pages with $2\leq r\leq n$,
according to Proposition \ref{spectral sequence trivial}. In turn, the associated graded on  the abutment is
$$
\gr H^i(Y,\O_Y) = \bigoplus_{r+s=i} H^s(\PS,\Lambda^r(R^1\varphi_*(\O_Y)).
$$
In our situation   $R^1\varphi_*(\O_Y)= \O_\PS(1)\oplus \shE_{-p} = F^*(\Omega^1(1))\oplus\O_\PS(1)$  by Proposition \ref{lie and image}, which 
shows that $\Lambda^r(R^1\varphi_*(\O_Y))=\shF_{r,0}\oplus\shF_{r-1,1}$, in the notation of Proposition \ref{exterior and frobenius}. 
Here we set $\shF_{-1,1}=0$ for convenience. Thus
\begin{equation}
\label{h^i formula}
h^i(\O_Y)=\sum_{s=0}^i  h^s(\shF_{i-s,0}) + \sum_{s=0}^i h^s(\shF_{i-s-1,1}).
\end{equation}
The summands for $s=0$ vanish by Proposition \ref{exterior and frobenius}. Furthermore, for $1\leq s\leq n-1$ we have 
$$
s\neq i-s\Rightarrow h^s(\shF_{i-s,0})=0\quadand s\neq i-s-1 \Rightarrow h^s(\shF_{i-s-1,1})=0.
$$
In the boundary case $i=n$, the 
last summands $h^n(\shF_{0,0})$ and $h^n(\shF_{-1,1})$ vanish,
the former by  \ref{exterior and frobenius}, the latter because   $\shF_{-1,1}=0$.
In turn,   the  sum \eqref{h^i formula} simplifies to 
$$
h^i(\O_Y) = \begin{cases}
h^j(\shF_{j,0})	 		&\text{if $i=2j$ is even;}\\
h^j(\shF_{j,1}) 		& \text{if $i=2j+1$ is odd.}
\end{cases}
$$
%Note that in both cases $j=\lfloor i/2\rfloor$.
The formula for the vector space dimensions now follows from Proposition \ref{exterior and frobenius}.
We also see that the filtration on the abutment $H^i(Y,\O_Y)$ has merely one step,
which gives the natural identification \eqref{natural identification} of groups.
\qed

\medskip 
Because of the relevance for the Picard scheme, we record:

\begin{corollary}
\mylabel{first cohomoloy groups}
For every $p>0$ and every $n\geq 0$ we have   $h^1(\O_Y)=n+1$, whereas
 $h^2(\O_Y)$ equals the number of lattice points $(l_0,\ldots,l_n)$
satisfying $l_0+\ldots+l_n=p$ and $0\leq l_0,\ldots,l_n\leq p-1$.  
\end{corollary} 

\proof
One easily checks that  terms   $H^i(\PS,\varphi_*(\O_Y))$, $i\geq 1$ and also the term $H^3(\PS, R^1\varphi_*(\O_Y))$
on the $E_2$-page for the Leray--Serre spectral sequence with respect to  $\varphi:Y\ra\PP^n$ vanish.
In turn, the formula for $h^1(\O_Y)$ and $h^2(\O_Y)$ of the Proposition hold regardless to the
assumption $p\geq n+1$.
In particular, $h^1(\O_Y)=n+1$, because the only lattice points in the polytope $P_{1,0}\subset\RR^{n+1}$ are the 
standard basis vectors.
\qed

\medskip
We now consider the case that gives varieties with $c_1=0$, which means 
$$
n=p-2,\quad \dim(Y)=2n+1=2p-3\quadand \omega_Y=\O_Y.
$$
With computer algebra \cite{Magma}, we computed the  cohomological invariants   for the first six   primes in the following table:
$$
\begin{array}{llll}
\toprule
p			& n 	& \dim(Y)	& h^0(\O_Y),\ldots,h^n(\O_Y)	\\
\midrule	
2			& 0		& 1	& 1	\\
\midrule
3			& 1		& 3	& 1, 2	\\
\midrule
5			& 3		& 7	& 1, 4, 52, 68\\
\midrule
7			& 5		& 11	& 1, 6, 786, 1251, 6891, 7872\\
\midrule
11			& 9		& 19	& 1, 10, 167950, 293830, 18480520, 25109950,\\
			& 		& 			& 251849140, 296659645, 859743835, 905642810\\
\midrule
13			& 11		& 23	& 1, 12, 2496132, 4457256, 825038490, 1149834280, 27258578260,\\
			& 		& 	& 33480335274, 223425722070, 250522227132, 616161367152,\\
			& 		& 	& 639330337978 \\
\bottomrule
\end{array}
$$
Note that the running time for $h^i(\O_Y)$ at the prime   $p=13$ and in degree $i=12$ was about three days.

%===========================================================
\section{Cohomology and weights}
\mylabel{Cohomology and weights}

The goal of this section is gain further control on the cohomology of
the Moret-Bailly family $Y=(A\times\PS)/H$, in particular for $H^2(Y,\O_Y)$ and $H^1(Y,\Theta_Y)$,
by using   automorphisms   of the superspecial abelian variety $A$ and their induced representations on
cohomology. This is best formulated with the machinery of weights, and requires some preparation.
Fix some integer $l\geq 1$, and consider 
$$
G=\mu_l=\GG_m[l]=\Spec\ZZ[T]/(T^l-1),
$$
viewed as a family of  finite diagonalizable group schemes
over the ground ring $R=\ZZ$.
Let $S$ be a scheme endowed with trivial $G$-action.
Recall that by \cite{SGA 3a}, Expos\'e I, Proposition 4.7.3,  
a \emph{$G$-linearization} for a quasicoherent sheaf $\shF$ on $S$ is   nothing but
a \emph{weight decomposition} $\shF=\bigoplus_w \shF_w$, where $w$ runs over the
character  group $\ZZ/l\ZZ=\uHom(G,\GG_m)$. Such characters are also known as \emph{weights}.
 A weight $w$ is called \emph{trivial} if $\shF_w=0$. We say  that the sheaf
$\shF$ is \emph{pure of weight $w_0$} if 
all other  weights $w\neq w_0$ are trivial.

For each ring $A$, the group elements
$\zeta\in G(A)\subset A^\times$   act   on $\shF_w\otimes_RA$ via   multiplication
by $\lambda=w(\zeta)$.
Note that  for every base-change $a:\Spec(k)\ra S$
for some field $k$  
containing a primitive $l$-th root of unity $\zeta$, 
the   weight decomposition for $\shF$ becomes   
  the \emph{eigenspace decomposition}  on the vector space $V=\shF(a)$ 
for the automorphism $\zeta:V\ra V$ with respect to the eigenvalues $\lambda=w(\zeta)$.
Also note that   $G$-linearizations for the sheaf $\shF$  correspond  to  $G$-actions on the finite $S$-scheme
stemming from the sheaf of dual numbers $\shA=\O_S\oplus\shF$, or   the vector $S$-scheme 
coming from   $\shA=\Sym^\bullet(\shF)$. 

Now assume we are over a ground field $k$ of  characteristic $p>0$. Since $\mu_l$ is already defined over $\FF_p$ we have
a canonical identification $\mu_l^\frob=\mu_l$, and the relative Frobenius map  
$\mu_l\ra\mu_l^\frob=\mu_l$ given by $ \zeta\mapsto \zeta^p$
induces multiplication by $p$ on the character group.
In turn, the Frobenius pullback $\shF^\frob$ has two canonical $G$-linearizations:
One stemming from $\mu_l^\frob=\mu_l$ with $(\shF^\frob)_w=(\shF_w)^\frob$,
the other via the Frobenius map, such that $ (\shF^\frob)_{pw}=(\shF_w)^\frob$.
It turns out that  the latter is more important for us.

In what follows we assume that $G=\mu_l$ acts on an abelian variety $A$ of dimension $g\geq 1$
such  that there is a $\mu_l$-equivariant principal polarization $A\ra P$, for the dual abelian variety $P=\Pic^0_{A/k}$
with the induced $\mu_l$-action. Note that the latter indeed can be achieved by passing to  $B=(A\oplus P)^{\oplus 4}$
and making a finite field extension (\cite{Koehler; Roessler 2003}, Lemma 3.2).
Furthermore, we assume that  the induced representation on $\Lie(A)$ is pure of weight $w_0\in\ZZ/l\ZZ$.
Using the equivariant bijections $H^1(A,\O_A)=\Lie(P)\ra \Lie(A)$ we immediately obtain:

\begin{lemma}
\mylabel{weight degree one}
The induced $\mu_l$-representation on $H^1(A,\O_A)$ is also pure of weight $w=w_0$.
\end{lemma}

Now suppose  additionally that $A$ is superspecial, and consider 
the Moret-Bailly family $Y=(A\times \PS)/H$
formed with    the homogeneous polynomials $Q_i=T_i$ of degree $d=1$, as in Section \ref{Cohomology and points}.
Since the induced action of $G=\mu_l$ on $\Lie(A)$ is pure, 
each copy $\alpha_p\subset A$ is normalized by $G$, and we get   induced actions
on the quotients $A/\alpha_p$. 

Consider the diagonal $G$-action on $X=A\times\PP^n$, with trivial action on the second factor.
This also can be seen as an action of the relative group scheme $G\times\PS=\mu_{l,\PS}$ that normalizes
the action of the family $H\subset X$ of height-one group schemes with $\Lie_{H/\PS}=\O_\PS(-1)$, and thus induces a  $G$-action
on our Moret-Bailly family $Y=X/H$. The structure morphism $\varphi:Y\ra\PS$ is equivariant,
with trivial $G$-action on the base. By Proposition \ref{lie sequence exact} we have a four-term exact sequence
$$
0\lra \O_\PS(-1)\lra \Lie(A)\otimes_k\O_\PS\lra \Lie_{Y/k}\lra \O_\PS(-p)\lra 0
$$
The cokernel for the inclusion on the left is 
the sheaf $\shE_1$,
and the above yields the  short exact sequence $0\ra\shE_1\ra \Lie_{Y/k}\ra\O_\PS(-p)\ra 0$
of sheaves with $G$-linearizations.   
We   saw in the proof for Proposition \ref{lie and image} that $\Ext^1(\O_\PS(-p),\shE_1)=0$.
By \cite{SGA 3a}, Expos\'e I, Proposition 4.7.4 we may choose a  splitting that respects  linearizations.
In turn, $\Lie_{Y/\PS} =\shE_1\oplus \O_\PS(-p)$
as  sheaves with $G$-linearization.

\begin{lemma}
\mylabel{weights in lie}
In the above setting, the summand  $\shE_1$ is pure of weight $w=w_0$, whereas $\O_\PS(-p)$ is pure of weight $pw_0$.
\end{lemma}

\proof
The sheaf  $\Lie(A)\otimes_k\O_\PS$ is pure of weight $w=w_0$, hence the same holds for the
subsheaf $\O_\PS(-1)$ and the quotient sheaf $\shE_1$. The relative Frobenius map 
$A\ra A^\frob$ becomes $G$-equivariant, provided that we take the induced $G$-action 
via the Frobenius map $\mu_l\ra\mu_l^\frob=\mu_l$.
In turn, the  $\O_\PS(-p)\subset\Lie(A^\frob)\otimes_k\O_\PS$ is pure of weight $w=pw_0$.
\qed

\medskip
We now make a similar analysis for the  higher direct images $R^i\varphi_*(\O_Y)$, which also 
come with   induced $G$-linearizations. 
To understand their weight decompositions, it suffices to treat the case $i=1$,
according to Proposition \ref{spectral sequence trivial}. Now we use the canonical projection
$Y=(A\times\PS)/H\ra A/A[F]\times\PS= A^\frob\times\PS$.
This map becomes $G$-equivariant if the right-hand side is endowed with the $\mu_l$-action
coming from the Frobenius map.

We already saw in Proposition \ref{lie and image} that we have an exact sequence of coherent sheaves $0\ra \O_\PS(1)\ra R^1\varphi_*(\O_Y)\ra \shE_{-p}\ra 0$,
and that all such extension split. As above, we see that there is a direct sum decomposition
$R^1\varphi_*(\O_Y) =\O_\PS(1)\oplus \shE_{-p}$
of sheaves with $G$-linearizations.

\begin{lemma}
\mylabel{weights in dual}
In the above setting, the summand 
 $\O_{\PP^n}(1)$ is pure of weight $w=pw_0$, whereas $\shE_{-p}$ is pure of weight $p^2w_0$.
\end{lemma}

\proof 
The Lie algebra $\liea=\Lie(A)$  is pure of weight $w=w_0$. According to Lemma \ref{weight degree one},
the same holds for $H^1(A,\O_A)$. In turn, the Frobenius pullback $H^1(A,\O_A)^\frob$
is pure of weight $w=pw_0$, since we use the action stemming from 
the Frobenius map $\mu_l\ra\mu_l^\frob$. We now proceed as for Lemma \ref{weights in lie}.
\qed

\medskip
The $G$-action on the Moret-Bailly family $Y$ 
induces $G$-representations on     Hodge groups   $H^s(Y,\Omega^r_Y)$ and the   tangent cohomology  
 $H^s(Y, \Theta_Y)$. Since the $G$-action respects the morphism $\varphi:Y\ra\PS$, we also
have induced representations on the groups
$H^s(Y,\varphi^*\Theta_\PS)$ and   $H^s(Y,\Theta_{Y/\PP^n})$. We now   compute the weights in some relevant cases:
 
\begin{proposition}
\mylabel{relevant weights}
Suppose we have $(p,n)\neq (2,1)$.  With respect to the induced $G$-representation, the following holds:
\begin{enumerate}
\item
The  group $H^2(Y,\O_Y)$ is pure of weight $w=p^2w_0$. 
\item
For  the   cohomology group $H^1(Y,\Theta_{Y/\PP^n})$ the only possible non-trivial weights are of the form  $w=mw_0$    with coefficient
$m\in\{p+1,p^2+1,p^2+p\}$.
\item 
For $H^1(Y,\varphi^*\Theta_Y)$ the only possible non-trivial weights are of the form $w=mw_0$ with $m\in\{p,p^2\}$.
\end{enumerate}
\end{proposition}

\proof
According to Proposition \ref{hodge numbers}, the Leray--Serre spectral sequence for the morphism $\varphi:Y\ra\PS$ induces
a canonical identification
$H^2(Y,\O_Y)=H^1(\PS,\shE_{-p})$. By  Lemma \ref{weights in dual}, the sheaf $\shE_{-p}$
is pure of weight $w=p^2w_0$, and assertion (i) follows.

Next we compute the weights in $H^1(Y,\Theta_{Y/\PP^n})$.
The projection formula for $\Theta_{Y/\PS}=\varphi^*\Lie_{Y/\PS}$ gives $R^i\varphi_*(\Theta_{Y/\PS})=\Lie_{Y/\PS}\otimes R^i\varphi_*(\O_Y)$,
and the Leray--Serre spectral sequence yields an exact sequence
$$
0\lra H^1(\PS, \Lie_{Y/\PS})\lra H^1(Y,\Theta_{Y/\PS}) \lra H^0(\PS, \Lie_{Y/\PS}\otimes R^1\varphi_*(\O_Y)).
$$
Recall from Proposition \ref{lie and image} that 
\begin{equation}
\label{splittings}
\Lie_{Y/\PS}=\shE_1\oplus\O_\PS(-p)\quadand R^1\varphi_*(\O_Y)=\O_\PS(1)\oplus \shE_{-p}.
\end{equation}
One easily computes $H^r(\PS,\shE_1)=0$  for all   $r\geq 1$ using the short exact sequence  
$0\ra\O_\PS(-1)\ra\O_\PS^{\oplus n+1}\ra\shE_1\ra 0$. Furthermore,
$H^1(\PS,\O_\PS(-p))$ is non-zero only if $n=1$ and $p=2$, which was  excluded.
In turn, it suffices to understand the possible weights in  $H^0(\PS, \shF)$ for the
sheaf $\shF= \Lie_{Y/\PS}\otimes R^1\varphi_*(\O_Y)$.  
For this we compute the weights occurring in each of  the summands
\begin{equation}
\label{summands}
\shE_1\otimes\O_\PS(1),\quad \O_\PS(-p)\otimes\O_\PS(1),\quad \shE_1\otimes\shE_{-p}  \quadand \O_\PS(-p)\otimes\shE_{-p} 
\end{equation}
inside $\shF$. According to Lemma \ref{weights in lie} and Lemma \ref{weights in dual}, these 
are pure of weight $w=mw_0$, where the coefficient $m$ is an integer of the form $1+p$ and $p+p$ and $1+p^2$
and $p+p^2$, respectively. The second case does not contribute, because the global sections for $\O_\PS(1-p)$ vanish.
Assertion (ii) follows.

It remains to understand the weights in  $H^1(Y,\varphi^*\Theta_\PS)$. 
Now the projection formula gives $R^i\varphi_*(\varphi^*\Theta_\PS)=\Theta_\PS\otimes R^i\varphi_*(\O_Y)$,
and the Leray--Serre spectral sequence yields an exact sequence
$$
0\lra H^1(\PS,\Theta_\PS)\lra H^1(Y,\varphi^*\Theta_\PS)\lra H^0(\PS,\Theta_\PS\otimes R^1\varphi_*(\O_Y)).
$$
From the Euler sequence $0\ra\O_\PS\ra\O_\PS^{\oplus n+1}(1)\ra\Theta_\PS\ra 0$ one easily infers that
the term on the left vanishes. As above, it suffices to understand the weights in  in  $H^0(\PS, \shF)$ for the coherent sheaf
$\shF=\Theta_\PS\otimes R^1\varphi_*(\O_Y)$. This is the sum of  
$
\Theta_\PS\otimes\O_\PS(1)\quadand \Theta_\PS\otimes\shE_{-p}
$.
Now the summands are pure of weight $w=mw_0$ with $m=0+p$ and $m=0+p^2$, which gives (iii).
\qed

\medskip
We now deduce a crucial fact that can be formulated without the machinery of weights: 

\begin{corollary}
\mylabel{sign involution}
Suppose that $p\neq 2$. Then the sign involution  on the superspecial abelian variety $A$
induces an action of the multiplicative group  $G=\{\pm 1\}$ on our Moret-Bailly family $Y=(A\times\PS)/H$,
and the induced representation on  the cohomology groups  $H^2(Y,\O_Y)$ and $H^1(Y,\varphi^*\Theta_Y)$ are multiplication by $\lambda=-1$,
whereas on  $H^1(Y,\Theta_{Y/\PP^n})$ it is multiplication by  $\lambda=1$.
\end{corollary}

\proof
We may assume that $k$ is algebraically closed,
and regard the abstract group $G=\{\pm 1\}$ as the  diagonalizable group scheme $G=\mu_l$ with $l=2$. 
Write $A=E_1\times\ldots\times E_g$ as a product of supersingular elliptic curve, and use
the canonical identification $E_i=\Pic_{E_i/k}$ to obtain an equivariant principal polarization.

The induced $G$-representation
on $\Lie(A)$ is multiplication by $\lambda=-1$. This is pure of weight $w_0=1$, seen as an element of the character group $\ZZ/2\ZZ$.
As discussed above, we get an induced action on $Y$, such that the quotient map $\epsilon: A\times\PS\ra Y$
is equivariant.
According to Proposition \ref{relevant weights}, the induced representation on $H^2(Y,\O_Y)$ is pure of weight $w=p^2 w_0=1$, because
$p\equiv 1$ modulo $l$. Similarly, $H^1(Y,\varphi^*\Theta_Y)$ has weight $w=1$.
In contrast, the only possible non-zero weights on $H^1(Y,\Theta_{Y/\PS})$ are $w=mw_0$,
with coefficient $m\equiv 0$  modulo $l$. In turn, $H^1(Y,\Theta_{Y/\PS})$ is pure of weight $w=0$.
\qed

%===========================================================
\section{Non-existence of formal liftings}
\mylabel{Non-existence formal}

We continue to study our Moret--Bailly families $Y=(A\times\PP^n)/H$ over the ground field $k$ of characteristic $p>0$,
formed with a superspecial abelian variety $A$ of dimension $g=n+1$, and the homogeneous polynomials $Q_i=T_i$ of degree $d=1$. 
The sign involution $a\mapsto -a$ on the abelian variety  $A$
induces an action of the cyclic group $G=\{\pm 1\}$ on the total space $Y$. For simplicity, the inclusion $G\subset\Aut(Y/k)$
is also called the \emph{sign involution}.
We now regard this $G$-action on $Y$  as additional structure,
and seek to understand its behavior under deformations.

Let $R$ be a complete local noetherian ring with residue field $R/\maxid_R=k$,
and $\foY\ra\Spf(R)$ be a proper flat morphism of formal schemes with closed fiber $Y=\foY\otimes_Rk$.
We say that this morphism is \emph{projectively algebraizable} if there is an invertible sheaf on $\foY$
whose restriction to $Y$ is ample. 
According to Grothendieck's Existence Theorem, the flat formal $R$-scheme $\foY$   then is the formal completion of some
flat projective $R$-scheme (\cite{EGA IIIa}, Theorem 5.4.5). 

\begin{proposition}
\mylabel{formal implies projective}
Suppose that $p\geq 3$. If the sign involution $G\subset\Aut(Y)$ extends to some $G\subset\Aut(\foY/R)$,
then $\foY\ra\Spf(R)$ is projectively algebraizable.
\end{proposition}

The proof for Proposition \ref{formal implies projective} requires a bit of  preparation and will be given below.
Let us first apply the result:
We say that the proper scheme $Y$ together with the sign involution  \emph{formally lifts  to characteristic zero}
if we may choose some $\foY\ra\Spf(R)$ and $G\subset\Aut(\foY/R)$  as above, where furthermore  the canonical map $\ZZ\ra R$ is injective.
Note that we may assume that $R$ is integral and one-dimensional,
by passing to the residue class ring for some suitable prime ideal.

Also note that  $A$, like any abelian variety,   projectively lifts to characteristic zero, but for $g\geq 2$ there are   formal liftings
$\foA\ra\Spf(R)$ to characteristic zero that are not algebraizable. The situation for   Moret-Bailly families is completely different.
From the above result, using Corollary \ref{dualizing sheaf} and Theorem \ref{no projective lift}, we immediately get:

\begin{theorem}
\mylabel{no formal lift}
Suppose that     $p\geq 3$ and  $p-1\leq g$.
Then the  Moret-Bailly family $Y=(A\times\PP^n)/H$   together with the sign involution   does not formally lift to characteristic zero.
\end{theorem}

We now prepare for the proof of Proposition \ref{formal implies projective}.
Recall that a \emph{$G$-linearized invertible sheaf} on $Y=(A\times\PP^n)/H$ is an invertible sheaf $\shL$,
together with a $G$-linearization, that is, a $G$-action on the line bundle
$L=\Spec(\Sym^\bullet(\shL^\vee))$ such that the structure morphism $L\ra Y$
is equivariant. Write $\Pic(Y,G)$ for the   abelian group of isomorphism classes for $G$-linearized invertible sheaves.
This can also  be viewed as the equivariant cohomology group $H^1(Y,G;\O_Y^\times)$, or the Picard group
of the quotient stack $[Y/G]$.

Now let $\foY\ra\Spf(R)$ be any  formal deformation of the Moret-Bailly family $Y$ over some complete local noetherian ring 
$R$ with residue field $R/\maxid_R=k$, and suppose that the sign involution  extends to some $G\subset\Aut(\foY/R)$.

\begin{proposition}
\mylabel{linearized sheaf extends}
The restriction map $\Pic(\foY,G)\ra\Pic(Y,G)$ is surjective.
\end{proposition}

\proof
Refining the descending chain $\maxid_R^j$, we get a descending chain
of ideals $\ideala_i$ with $\length(\ideala_i/\ideala_{i+1})=1$, such that $R=\invlim R_i$,
where $R_i=R/\ideala_i$. Setting   $Y_i=\foY\otimes_RR_i$, we get an
increasing sequence $Y=Y_0\subset Y_1\subset \ldots$, where each $Y_i$ is a proper flat $R_i$-scheme,
and the comparison maps $Y_{i-1}\ra Y_i\otimes_{R_i}R_{i-1}$ are isomorphisms.
In fact, the formal scheme $\foY$ is nothing but  the resulting inverse system $(Y_i)_{i\geq 0}$ of $R$-schemes.

Let $\shL_0$ be a $G$-linearized invertible sheaf on $Y_0=Y$.
We show by induction on $i\geq 0$ that it extends to $Y_i$. This is obvious for $i=0$.
Suppose now that  we have a linearized extension $\shL_i$ on $Y_i$.
The ideal sheaf $\shI$ for the closed embedding $Y_i\subset Y_{i+1}$ has square zero, and is isomorphic to $\O_Y$ 
by flatness. In turn, we have a short exact sequence $0\ra\O_Y\ra\O_{Y_{i+1}}^\times\ra\O_{Y_i}^\times\ra 1$, which gives 
an exact sequence
$$
0\lra H^1(Y,\O_Y)\lra \Pic(Y_{i+1})\lra \Pic(Y_i)\stackrel{\partial}{\lra} H^2(Y,\O_Y).
$$
By naturality of cohomology, this  sequence is equivariant with respect to the induced $G$-actions.
The isomorphism class of $\shL_i$ is $G$-invariant, so Corollary \ref{sign involution} gives
$\partial(\shL_i)=-\partial(\shL_i)$. With $p\neq 2$ we infer that the obstruction $\partial(\shL_i)$ vanishes,
hence $\shL_i$ extends to some invertible sheaf on $Y_{i+1}$. 

It remains to choose an extension that admits a linearization.
To achieve this we use the results form \cite{Schroeer; Takayama 2018}, discussed in more details in the   proof for Proposition \ref{equivalences for lifting} below.
Let $L$ be the set of isomorphism classes for pairs $(\shL',\varphi')$, where $\shL'$
is invertible on  $Y_{i+1}$ and $\varphi':\shL_i\ra\shL'|Y_i$ is an isomorphism,
and define $T=H^1(Y,\O_Y)$. Then $L$ carries the structure of an $T$-torsor with group of operators $G$, giving
a cohomology class $[L]\in H^1(G,T)$. This cohomology group vanishes because $|G|=2$ is relatively prime to the characteristic $p\geq 3$.
In turn, we can choose a   lifting $\shL_{i+1}$ whose isomorphism class is $G$-fixed.
The group of automorphisms of $\shL_{i+1}$ restricting to the identify on $\shL_i$
is the vector space $V=\Hom(\shL_{i+1},\O_Y)=\Hom(\shL_0,\O_Y)$.
According to \cite{Schroeer; Takayama 2018}, Theorem 1.2 the   obstruction against the existence of a $G$-linearization lies in the cohomology
group $H^2(G,\ideala_i/\ideala_{i+1}\otimes_kV)$, which again vanishes.
Summing up, the linearized sheaf $\shL_i$ extends to $Y_{i+1}$.
\qed

\medskip
\emph{Proof for Proposition \ref{formal implies projective}.}
Suppose  there is a    formal lifting $\foY\ra\Spf(R)$ to characteristic zero,
such that the sign involution extends to some $G\subset\Aut(\foY/R)$.
Choose a very ample invertible sheaf $\shL$ on $Y$.  
Replacing $\shL$ by $\shL\otimes\sigma^*(\shL)$, where $\sigma\in G$ is the generator,
we may assume that the isomorphism class is $G$-invariant.
According to \cite{Grothendieck 1957}, Theorem 5.2.1 there is a  spectral sequence
$$
H^r(G,H^s(Y,\O_Y^\times))\Longrightarrow H^{r+s}(Y,G,\O_Y^\times).
$$
The resulting five term exact sequence shows that  the obstruction for the existence of  a $G$-linearization on $\shL$ lies in the elementary abelian group
$H^2(G,k^\times)=k^\times/k^{2\times}$. Thus $ \O_Y(1)=\shL^{\otimes 2}$   can be endowed with
a $G$-linearization, and is very ample.

By Proposition \ref{linearized sheaf extends} there is some
invertible sheaf $\O_\foY(1)$ restricting to $\O_Y(1)$.  
According to Grothendieck's Existence Theorem,  $\foY$ is the formal completion of some
flat projective $R$-scheme (\cite{EGA IIIa}, Theorem 5.4.5). In turn, our Moret--Bailly family $Y$ admits a  projective lift
to characteristic zero.
\qed

\medskip
Recall that the Moret--Bailly family $Y=(A\times\PP^n)/H$ comes with an induced morphism $\varphi:Y\ra\PP^n$.
It turns out that extending the sign involution is essentially the same as extending this morphism.
Let $\foY\ra\Spf(R)$ be any   formal deformation of the Moret-Bailly family $Y$ over some complete local noetherian ring 
$R$ with residue field $R/\maxid_R=k$

\begin{proposition}
\mylabel{equivalences for lifting}
Suppose $p\geq 3$. Then the following two conditions are equivalent:
\begin{enumerate}
\item
The     morphism $\varphi:Y\ra\PP^n$ extends to a morphism $\foY\ra\PP^n_R$.
\item
The sign involution $G\subset \Aut(Y/k)$ extends to  an inclusion $G\subset \Aut(\foY/R)$.
\end{enumerate}
\end{proposition}

\proof
We start with the implication (i)$\Rightarrow$(ii).
Suppose   that  $\varphi:Y\ra \PP^n$ extends to a morphism, which we likewise 
call $\varphi:\foY\ra\PP^n_R$. In order to extend the action of $G=\{\pm1\}$,
we use Rim's results about equivariant structures on versal deformations \cite{Rim 1980}, in the form obtained 
in \cite{Schroeer; Takayama 2018}.  Choose a Cohen ring $\Lambda$ with residue field $k$. Recall that this a complete
discrete valuation ring with maximal ideal $\maxid_\Lambda=p\Lambda$.
Let $(\Art_\Lambda)$ be the category of local Artin rings $A$ over $\Lambda$ with residue field $k$,
and $\shF\ra(\Art_\Lambda)^\op$ be the category fibered in groupoids whose fiber category over $A$
comprises the  pairs $\xi=(X,\alpha)$, where $X$ is a proper flat  scheme over $\PP^n_{R}$, and $\alpha:Y\ra X\otimes_Rk$
is an isomorphism. This is a \emph{deformation category} in the sense of Talpo and Vistoli \cite{Talpo; Vistoli 2013}.
To conform with  \cite{Schroeer; Takayama 2018}, the symbol $A$, which otherwise denotes our superspecial  abelian variety,
is here also used for local Artin rings; this should not cause any confusion.

Using the notation from the proof of Proposition \ref{linearized sheaf extends}, we write $R=\invlim R_i$.
We now show  by induction on $i\geq 0$ that there are     compatible inclusions $G\subset\Aut(Y_i/\PP^n_{R_i})$.
For $i=0$ we choose the given action on $Y_0=Y$.
Now suppose that we already have the action on $Y_i$.
Using the notation from \cite{Schroeer; Takayama 2018},  set 
$A=R_i$ and $A'=R_{i+1}$. Write $\xi \in \shF(A)$ for the given family $Y_i=\foY\otimes_RA$ over $\PP^n_A$,
and $\Lif(\xi,A')$ for the set of isomorphism classes $[f]$ of all cartesian morphisms $f:\xi\ra \xi'$ over the inclusion $\Spec(A)\subset\Spec(A')$,
in the fibered category $\shF$.
This set comes with a  $G$-action, namely $\sigma \cdot[f]=[f\circ\sigma^{-1}]$.
The family  $Y_{i+1}=\foY\otimes_RA'$ yields an element $\xi_{i+1}$ in   $L=\Lif(\xi,A')$, which thus is non-empty.
Furthermore,  it carries the structure of  a torsor for the \emph{tangent space} 
$T=H^1(Y,\Theta_{Y/\PP^n})$ for the deformation category. It is also endowed with $G$ as  group of operators,
meaning that the action
\begin{equation}
\label{torsor}
L\times T\lra L,\quad ([f],t)\longmapsto [f]\cdot t
\end{equation}
satisfies ${}^{\sigma}([f]\cdot t)={}^\sigma[f]\cdot {}^\sigma t$.  

Now write   $f:\xi\ra \xi'$ for the lifting corresponding
to $Y_{i+1}$, and let $\sigma\in G$ be the generator.
Since $\sigma:Y_i\ra Y_i$ is compatible with the morphism $\varphi:Y_i\ra\PP^n_{R_i}$,
the new class $[f\circ \sigma^{-1}]$  differs from the old class by some element $t\in T$, in other word
${}^\sigma[f] = [f]\cdot t$. Applying the involution to this equation, and using that it acts trivially on $T$ by our weight
computation in Proposition \ref{sign involution}, we get
$$
[f] = {}^\sigma([f]\cdot t) = {}^\sigma [f] \cdot {}^\sigma t = {}^\sigma [f] \cdot t = [f]\cdot 2t.
$$
This ensures $2t=0$, and hence $t=0$ since we are in characteristic $p\neq 2$.
In other words, $[f]\in L$ is   $G$-fixed.

Consequently, the object $\xi'\in\shF(A')$ corresponding to   $Y_{i+1}=\foY\otimes_RA'$ has $G$-fixed isomorphism class.
By loc.\ cit., Theorem 3.3 the obstruction for extending the $G$-action from $Y_i$  to $Y_{i+1}$
lies in $H^2(G, \ideala_i/\ideala_{i+1}\otimes_k\Aut_{\xi_0}(\xi_{k[\epsilon]}))$.
This group   vanishes  because $p\neq 2$, and we conclude that the  action on $Y_i$ extends to an  action on $Y_{i+1}$ over $\PP^n_{R_{i+1}}$.
 
We now establish the reverse implication (ii)$\Rightarrow$(i). Suppose that we have $G\subset\Aut(\foY/R)$.
By induction on $i\geq0$, we show that $Y_i\ra\Spec(R_i)$ factors over $\PP^n_{R_i}$
in a compatible way, such that $G$ acts over $\PP^n_{R_i}$.
 
Suppose we already have a factorization $Y_i\ra\PP^n_{R_i}$. This morphism must be flat, by the fiber-wise   criterion for flatness 
(\cite{EGA IVc}, Theorem 11.3.10).
The obstruction to extend $Y_i$ over $\PP^n_{R_{i+1}}$ is an element $\alpha\in H^2(Y,\Theta_{Y/\PP^n})$.
The short exact sequence $0\ra\Theta_{Y/\PP^n}\ra \Theta_Y\ra \varphi^*\Theta_{\PP^n}\ra 0$ induces an exact sequence
$$
H^1(Y,\varphi^*\Theta_{\PP^n})\lra H^2(Y,\Theta_{Y/\PP^n})\lra H^2(Y,\Theta_Y).
$$
Since $Y_{i+1}$ extends $Y_i$ over $R_{i+1}$, the image of $\alpha$ in $H^2(Y,\Theta_Y)$ vanishes, hence our obstruction belongs to the image of $H^1(Y,\varphi^*\Theta_{\PP^n})$.
Let $\sigma\in G$ be the generator. By the weight computation in Corollary \ref{sign involution}, we must have  $\sigma^*(\alpha)=-\alpha$.
On the other hand, the functoriality of the obstruction gives $\sigma^*(\alpha)=\alpha$. Using $p\neq 2$ we conclude that the obstruction
vanishes, thus we have an extension $Y_{i+1}\ra \PP^n_{R_{i+1}}$.

It remains to check that the $G$-action on $Y_{i+1}$ is over $\PP^n_{R_{i+1}}$. 
In light of the implication (i)$\Rightarrow$(ii), it suffices to check that there is at most extension of $G\subset\Aut(Y_i/R_i)$
to $G\subset\Aut(Y_{i+1}/R_{i+1})$. 
As explained in \cite{Schroeer; Takayama 2018}, page 412 any two such extension differ by a homomorphism of $G$ into the vector space
of automorphisms of $Y_{i+1}$ restricting to the identity on $Y_i$. Using $p\neq 2$, we see that such homomorphisms are trivial. 
\qed

%===========================================================
\section{Non-existence of liftings to Witt vectors}
\mylabel{Non-existence witt}

Fix integers $n,d\geq 1$  and form the Moret--Bailly family $Y$  with some superspecial abelian variety $A$
of dimension $g=n+1$ over a field $k$
of characteristic $p>0$, as in Section \ref{Moret-Bailly}.
Recall that it comes with two morphisms $\varphi:Y\ra\PS$ and $\psi:Y\ra A^{(p)}$. Finally, 
recall that $\dim(Y)=2n+1$ and that $\omega_Y=\varphi^*\O_\PS(m)$, where $m=d(p-1) - (n+1)$. 
Here $\omega_Y=\det(\Omega^1_{Y/k})$ is the dualizing sheaf of $Y$.

We suspect that regardless of
the values of  $n,d$ the scheme $Y$ does not lift to characteristic zero,
and perhaps also not to the ring $W_2(k)$ of Witt vectors of length two.
Regarding $Y=Y_{n,d,p}$ also in dependence of the prime $p=\operatorname{char}(k)$, we will show:

\begin{theorem}
\mylabel{no lift to witt}
Fix $n\geq 2$ and $d\geq 1$. Suppose that  $n\not\equiv 2$ modulo $4$ and   that the ground field $k$ is perfect. Then
the Moret--Bailly family $Y$   
does not lift over  the ring $W_2(k)$, for    almost all primes $p>0$.
\end{theorem}

This is an application  of the  Deligne--Illusie result  on    Kodaira--Nakano--Akizuki Vanishing in positive characteristics,
together with a   Hirzebruch--Riemann--Roch computation    of certain Euler characteristics based on our knowledge
of the tangent sheaf $\Theta_{Y}$. The proof requires some preparations and is given
 towards the end of the section.

Recall that the Chow ring with rational coefficients of the projective space is a truncated polynomial ring
$\CH^\bullet(\PP^n)=\QQ[h]/(h^{n+1})$, where $h\in\CH^1(\PP^n)$ is the class of the invertible sheaf $\O_{\PP^n}(1)$.
Consider the formal power series
$$
Q(x) = \left(\frac{-dx}{1-e^{dx}}\right)^{-1} \left(\frac{-dpx}{1-e^{dpx}}\right)  \left(\frac{x}{1-e^{-x}}\right)^{n+1}\in\QQ[[x]].
$$
It induces an element $Q(h)\in \CH^\bullet(\PP^n)$, which  appears in the following computation:

\begin{proposition}
\mylabel{hirzebruch-riemann-roch}
Let $\shE$ and $\shF$ be coherent sheaves on $A^{(p)}$ and $\PS$, respectively.
Then the coherent sheaf $\shM=\psi^*(\shE)\otimes\varphi^*(\shF)$ on the Moret--Bailly family $Y$ has Euler characteristic
 $$
\chi(\shM) = p^{g-1} \chi(\shE) \int_{\PP^n}  \ch(\shF)Q(h).
$$
\end{proposition}

\proof
We have  $\chi(\shM) = \int_Y\ch(\shM)\td(\Theta_{Y}) $ according to   Hirzebruch--Riemann--Roch,
where $\ch(\shM)$ is the Chern character of $\shM$,
 and $\td(\Theta_{Y})$ is the Todd class of the tangent sheaf (\cite{Fulton 1998}, Corollary 15.2.2). Recall that   Todd classes are
multiplicative in short exact sequences, hence 
$\td(\Theta_Y)=\td(\Theta_{Y/\PS})\td(\varphi^*\Theta_\PS)$.
Using the decomposition $\Theta_{Y/\PP^n} = \varphi^*(\shE_d\oplus \O_{\PP^n}(-dp))$ from Proposition \ref{lie and image} and the exact sequences
$$
0\ra \O_\PS(-d)\ra\O_\PS^{\oplus n+1}\ra\shE_d\ra 0 \quadand 0\ra \O_{\PP^n}\ra\O_{\PP^n}^{\oplus n+1}(1)\ra \Theta_{\PP^n}\ra 0,
$$   
we infer that the Todd class $\td(\Theta_Y)\in\CH^\bullet(Y)$ is the pullback  of
$$
Q(h)=\td(\shL^{\otimes-d})^{-1} \td(\shL^{\otimes -dp})\td(\shL)^{n+1}\in\CH^\bullet(\PS)
$$
with respect to $\varphi:Y\ra\PP^n$. Here  $\shL=\O_{\PP^n}(1)$, with first Chern class $h\in \CH^2(\PP^n)$ and 
Todd class  $\td(\shL) = h/(1-e^{-h})$.

The morphism $\varphi:Y\ra\PP^n$ factors over the morphism $\psi\times\varphi:Y\ra A^{(p)}\times\PP^n$ to the product, which is locally free
of rank $p^{g-1}$, and we have $\shM = (\psi\times\varphi)^*(\pr_1^*\shE\otimes\pr_2^*\shF)$. The projection formula thus  gives
$$
\int_Y\ch(\shM)\psi^*(Q(h)) = p^{g-1}\int_{A^{(p)}\times\PP^n} \ch(\pr_1^*\shE\otimes \pr_2^*\shF)\pr_2^*(Q(h)).
$$ 
Since  Chern characters are natural, and multiplicative in tensor products,   the above can be written as
$$
p^{g-1}\int_{A^{(p)}\times\PP^n} \pr_1^*(\ch(\shE)) \cdot \pr_2^*(\ch(\shF)Q(h)).
$$
By the functoriality of proper push-forwards and the projection formula,  this coincides with 
$$
p^{g-1}\int_{\PP^n}  \pr_{2*}(\pr_1^*(\ch(\shE)))\cdot \ch(\shF) Q(h).
$$
We now use that    the Todd class 
for the tangent sheaf of any abelian variety is the unit element in the Chow group.
Together with  the compatibility of proper push-forwards with flat pull-backs (\cite{Fulton 1998}, Proposition 1.7) applied to  the  the diagram
$$
\begin{CD}
A^{(p)}\times\PS	@>\pr_1>>	A^{(p)}\\
@V\pr_2VV			@VVfV\\
\PS		@>>g>	\Spec(k),
\end{CD}
$$
we see  
$\pr_{2*}(\pr_1^*(\ch(\shE) )) = g^*f_*(\ch(\shE) \td(\Theta_{A^{(p)}}))=\chi(\shE)$
 in the Chow group $\CH^\bullet(\PS)$.
This yields the desired formula for the Euler characteristic of $\shM$.
\qed

\medskip
Suppose now that $\shE$ is an ample invertible sheaf on $A^{(p)}$, and let $s\geq 1$ be an integer.
Since $\psi\times\varphi:Y\ra A^{(p)}\times\PS$ is finite, the invertible sheaf 
$\shL=\psi^*(\shE)\otimes\varphi^*(\O_{\PP^n}(s))$  on the Moret-Bailly family $Y$ is ample.

\begin{proposition}
\mylabel{negative euler characteristic}
In the above setting, fix the integers $n\geq 2$, $d\geq 1$ and some $s>(n+d+1)/2$. Furthermore  assume   $n\not\equiv 2$ modulo $4$. Then  the invertible sheaf  
$\shM=\shL\otimes\omega_Y$ has Euler characteristic $\chi(\shM)<0$ for almost all primes $p>0$.
\end{proposition}

\proof
 We have $\shM=\psi^*(\shE)\otimes\varphi^*(\shF)$
with $\shF=\O_\PS(d(p-1) - (n+1) + s)$, by  Proposition \ref{dualizing sheaf}. Hence  Proposition \ref{hirzebruch-riemann-roch} gives
$$
\chi(\shM) = p^{g-1} \chi(\shE) \int_\PS Q(h) e^{(dp-d - n-1 +s)h}.
$$
The factor $\chi(\shE)$ is strictly positive,   by Hirzebruch--Riemann--Roch  or \cite{Mumford 1970}, Section 16. Consequently, to 
determine the sign of $\chi(\shM)$ 
we merely have to understand the sign of the integrand in   top degree $n$, that is,  the sign of the
coefficient at $x^n$ in the formal power series
$$
\left(\frac{-dx}{1-e^{dx}}\right)^{-1} \left(\frac{-dpx}{1-e^{dpx}}\right)  \left(\frac{x}{1-e^{-x}}\right)^{n+1} e^{(dp-d - n-1 +s)x} .
$$
This series can be rewritten as
\begin{equation}
\label{formal power series}
\left(\frac{-dx}{1-e^{dx}}\right)^{-1} \left(\frac{x}{1-e^{-x}}\right)^{n+1}e^{-(d + n+1-s)h} \cdot  \left(\frac{dpx}{1-e^{-dpx}}\right),
\end{equation}
where the prime $p>0$ only enters through the last factor. We temporarily regard   $p$ as another indeterminate, and the above expression 
as element $\sum\lambda_i(p)x^i$ in the ring $\QQ[p][[x]]$. It appears difficult to find a closed formula
for the polynomial $\lambda_n(p)$. However, it is possible to determine the  sign of the leading
coefficient:

Recall that the formal power series
$x/(1-e^{-x}) = 1 + x/2 + \sum_{i=2}^\infty B_i x^i/i!$
can be expressed in terms of  the Bernoulli numbers $B_i\in\QQ$. We follow the convention that $B_1=1/2$. 
Since $n\geq 2$, the last factor in \eqref{formal power series} becomes
$$
1+ \frac{dp}{2}x +  \ldots + \frac{(dp)^{n-1}B_{n-1}}{(n-1)!}x^{n-1} + \frac{(dp)^nB_n}{n!}x^n   +\dots,
$$
while the product of the first three factors in \eqref{formal power series} starts with
$$
\left(1-\frac{dx}{2} + \ldots\right)^{-1} \left(1+\frac{x}{2}+\ldots\right)^{n+1} \left(1-\frac{(d+n+1-s)x}{1!}+\ldots \right).
$$
This equals  $1+\frac{2s-d-n-1}{2}x$ modulo $x^2$.
The coefficient  of $x^n$ in the power series \eqref{formal power series} thus takes the form
\begin{equation}
\label{polynomial in p}
\lambda_n(p) = \frac{B_n(dp)^n}{n!} + \frac{(2s-d-n-1)B_{n-1}(dp)^{n-1}}{2(n-1)!} + \text{lower degree}.
\end{equation}
By assumption, $2s-d-n-1>0$. The signs of the Bernoulli numbers are given by 
$$
\operatorname{sign}(B_{2j})=(-1)^{j-1}\quadand B_{2j+1}=0,
$$
for every $j\geq 1$. Note that the former  can be deduced from   the formula $B_{2j}=\frac{ (-1)^{j-1}2(2j)!} {(2\pi)^{2j}}\zeta(2j)$
involving the Riemann Zeta function $\zeta(z)$. 
It follows that  the degree of  $\lambda_n(p)$ depends on the parity of $n\geq 2$.
However,  for both $n=2j$ even and   $n=2j+1$ odd, the sign of the leading term in \eqref{polynomial in p} equals $(-1)^{j-1}$.
In other words, for $n\not\equiv 2$ modulo $4$ this sign is always negative. 
Regarding the symbol $p$
again as a prime number, we conclude that  $\lambda_n(p)<0$ for all sufficiently large primes $p\gg 0$,
and then $\chi(\shM)<0$. 
\qed

\medskip
\emph{Proof for Theorem \ref{no lift to witt}.}
Seeking a contradiction, we assume that  there are infinitely many primes $p>0$ 
for which the Moret--Bailly family $Y$ lifts to the ring of truncated Witt vectors $W_2(k)$,
with fixed integers $n\geq 2$ and $d\geq 1$.
Choose some $s> (n+d+1)/2$.
By Proposition \ref{negative euler characteristic} we find some sufficiently large $p\geq\dim(Y) $ so that there is an  ample invertible sheaf 
$\shL$ on $Y$ with 
$\chi(\shL\otimes\omega_Y)<0$. Consequently
 $h^i(\shL\otimes\omega_Y)>0$ for some odd degree $i$.
On the other hand, the  Deligne--Illusie form of Kodaira--Akizuki--Nakano Vanishing 
(\cite{Deligne; Illusie 1987}, Corollary 2.8) ensures that  $h^j(\shL\otimes\omega_Y)=0$ for all $j>0$,
contradiction.
\qed

\medskip
Of course, for fixed values of $n\geq 1$ one may compute   the coefficient $\lambda_n(p)$  in the formal 
power series \eqref{formal power series}, either  by hand or with computer algebra, as a polynomial in $d,p,s$. For $n=1$
we have $\lambda_1(p)= d(p-1)/2+s-1$, which is always positive.
For  $n=3$ and $s=1$ we obtain
$$
\lambda_3(p) = -\frac{d^3+2d^2}{24} p^2  + \frac{d^3+3d^2+2d}{12} p - \frac{d^3+2d^2+2d}{24}.
$$
The zero of the derivative $\lambda_3'(p)$ is the rational number $r=(d^3+3d^2+2d)/(d^3+2d^2)$,
which for all $d\geq 1$ satisfies $r<2$. Arguing as above, we obtain:
 
\begin{proposition}
\mylabel{no lift special case}
For $n=3$ and  $d\geq 1$, the Moret--Bailly family $Y$ over a perfect field of characteristic $p\geq 7=\dim(Y)$
does not lift to the ring of   Witt vectors
$W_2(k)$ of length two.
\end{proposition}

%===========================================================

\end{document}